\documentclass[a4paper,10pt,oneside,onecolumn,number,preprint,centertitle]{article}

\usepackage{amssymb, amsfonts, mathtools, systeme}
\usepackage[caption=false]{subfig}
\usepackage{xcolor}
\usepackage{multirow}
\usepackage{geometry}
\usepackage[colorlinks=true,linkcolor=blue,citecolor=red]{hyperref}
\AtBeginDocument{\hypersetup{citecolor=red}}
\usepackage{cleveref}
\crefname{figure}{Fig.}{Figs.}
\Crefname{figure}{Figure}{Figures}
\crefname{section}{Sec.}{Secs}

\newcommand{\x}{\boldsymbol{x}}
\newcommand{\pa}{\partial\Omega}
\newcommand{\Mp}{\mathcal{M}_p}
\def\ctanh{\mathrm{ctanh}}

\begin{document}

\title{A numerical study of the Dirichlet-to-Neumann operator in planar domains}
\author{Adrien Chaigneau$^{1,}$\thanks{adrien.chaigneau@polytechnique.edu}
	\and Denis S. Grebenkov$^{1,}$\thanks{denis.grebenkov@polytechnique.edu (Corresponding author)}}
\date{$^1$\textit{Laboratoire de Physique de la Matière Condensée, CNRS, Ecole Polytechnique, Institut Polytechnique de Paris, 91120, Palaiseau, France} \\[2ex]
\today}
%
\maketitle

\begin{abstract}
We numerically investigate the generalized Steklov problem for the modified Helmholtz equation and focus on the relation between its spectrum and the geometric structure of the domain. We address three distinct aspects: (i) the asymptotic behavior of eigenvalues for polygonal domains; (ii) the dependence of the integrals of eigenfunctions on the domain symmetries; and (iii) the localization and exponential decay of Steklov eigenfunctions away from the boundary for smooth shapes and in the presence of corners. For this purpose, we implemented two complementary numerical methods to compute the eigenvalues and eigenfunctions of the associated Dirichlet-to-Neumann operator for various simply-connected planar domains. We also discuss applications of the obtained results in the theory of diffusion-controlled reactions and formulate several conjectures with relevance in spectral geometry.
\end{abstract}
\noindent
\section{Introduction}\label{sec:intro}

The Dirichlet-to-Neumann operator \cite{taylor1996partial, agranovich2006mixed, girouard2017spectral, levitin2023topics} plays a prominent role in applied mathematics, physics, engineering and medicine. One of its most common applications is related to medical imaging and electrical impedance tomography \cite{cheney1999electrical,borcea2002electrical}, in which the electric conductivity in the bulk has to be determined from electric measurements on the boundary and allows for instance lung function assessment. A similar technique was used in geophysics for imaging sub-surface structures \cite{zhdanov1994geoelectrical,uhlmann2014inverse}. A recently developed theoretical description of diffusion-controlled reactions relies on the eigenbasis of the Dirichlet-to-Neumann operator to decompose the underlying propagators \cite{grebenkov2019spectral,grebenkov2020paradigm, grebenkov2020surface}. Althought the spectral properties of the Dirichlet-to-Neumann operator have been intensively studied over the past century, there are still many open questions and unsolved problems that explains a rapid development of this topic during the last years \cite{hislop2001spectral,arendt2011dirichlet,arendt2015dirichlet,polterovich2019nodal,girouard2017spectral,galkowski2019pointwise,daude2021exponential,helffer2022semi,levitin2022sloshing,girouard2022dirichlet,colbois2023some}.

In this study, we focus on simply-connected planar bounded domains $\Omega\in \mathbb{R}^2$ with a Lipschitz boundary $\pa$. The Dirichlet-to-Neumann operator $\Mp$ associates to a function $f$ on the boundary $\pa$ another function on that boundary:

\begin{equation}
\begin{aligned}
\Mp: H^{1 / 2}(\pa) &\rightarrow H^{-1 / 2}(\pa) \\
f &\mapsto \left.\left(\partial_n u\right)\right|_{\pa},
\end{aligned}
\end{equation}
where $\partial_n$ is the normal derivative oriented outward the domain and $u(\x)$ is the solution of the boundary value problem,

\begin{equation}\label{prob:helm}
\left\{ \begin{array}{r l l}
(p-\Delta) u(\x)& = 0 & \quad (\x \in \Omega), \\
u(\x)& = f(\x) &\quad  (\x \in \pa),
\end{array}  \right.
\end{equation}
in the Sobolev space 
\begin{equation}
\mathcal{H}^1(\Omega) = \{ u \in L^2(\Omega) \mid \partial_{x}u \in L^2(\Omega),  \partial_y u \in L^2(\Omega)\},
\end{equation}
where $\Delta=\partial_x^2+\partial_y^2$ is the Laplace operator, $p\in\mathbb{R}$ is a fixed parameter and $L^2(\Omega)$ is the space of measurable and square-integrable functions on $\Omega$. The functional space $H^{1 / 2}(\partial \Omega)$ is the trace of $H^{1}(\Omega)$:
\begin{equation}
\begin{aligned}
H^{1 / 2}(\partial \Omega):=&\operatorname{tr}\left(H^1(\Omega)\right) \\
:=&\left\{v \in L^2(\partial \Omega) \mid \exists u \in H^1(\Omega): \operatorname{tr}(u)= u|_{\pa} = v\right\},
\end{aligned}
\end{equation}
and the space $H^{-1/2}(\pa)$ is the dual of $H^{1/2}(\pa)$ (see details on functional spaces in \cite{brezis2011functional,levitin2023topics}).
In the context of diffusion-controlled reactions, the function $f$ can be thought of as a source of molecules on the boundary $\pa$, so that $\mathcal{M}_pf$ gives their flux density on that boundary.
Throughout the paper, we focus on $p\geq0$ so that $\Mp$ is a selfadjoint operator that has a discrete spectrum \cite{levitin2023topics}, with the eigenvalues $\mu_k^{(p)}$ and eigenfunctions $v_k^{(p)}$ satisfying 

\begin{equation}\label{eq:eigDtN}
\Mp v_k^{(p)} = \mu_k^{(p)} v_k^{(p)} \quad (k=0,1,2,\ldots).
\end{equation}
The eigenvalues are nonnegative and the eigenfunctions form an orthonormal basis of $L^2(\pa)$.
We enumerate the eigenvalues in increasing order 
\begin{equation}\label{eq:ord_eigv}
\mu_0^{(p)}\leq\mu_1^{(p)}\leq...\nearrow \infty.
\end{equation} 
The spectrum of $\Mp$ is closely related to the spectrum of the (generalized) Steklov problem \cite{stekloff1902problemes,kuznetsov2014legacy}:
\begin{equation}\label{prob:VkSteklov}
\left\{ \begin{array}{r l l}
(p-\Delta) V_k^{(p)} & = 0 & \quad (\x \in \Omega), \\
\partial_n V_k^{(p)}& = \mu_k^{(p)} V_k^{(p)} &\quad  (\x \in \pa),
\end{array}  \right.
\end{equation}
where $V_k^{(p)}$ are the Steklov eigenfunctions.
One sees that $v_k^{(p)}$ is the restriction of $V_k^{(p)}$ on $\pa$, whereas $V_k^{(p)}$ can be obtained as the unique extension of $v_k^{(p)}$ into $\Omega$:
\begin{equation}\label{prob:Vkextension}
\left\{ \begin{array}{r l l}
(p-\Delta) V_k^{(p)} & = 0 & \quad (\x \in \Omega), \\
V_k^{(p)}& = v_k^{(p)} &\quad  (\x \in \pa).	
\end{array}  \right.
\end{equation}

In this paper, we address three questions about the spectral properties of the Dirichlet-to-Neumann operator. First, we investigate the asymptotic behavior of the eigenvalues $\mu_{k}^{(p)}$ as $p\to\infty$. According to \cite{girouard2022dirichlet}, one has
\begin{equation}\label{eq:asymp-mu_circle}
\mu_{k}^{(p)} \simeq \sqrt{p} \quad(p\gg1),
\end{equation}
for all $k$ and all bounded domains with a smooth $\mathcal{C}^1$ boundary. The symbol $\simeq$ denotes the asymptotic behavior of $\mu_k^{(p)}$ when $p$ goes to infinity; however, it also
emphasizes that the left-hand side is close to the right-hand side when $p$ is large enough. In turn, in the presence of corners one can expect 
\begin{equation}\label{eq:asymp-ck}
\mu_{k}^{(p)} \simeq  c_k\sqrt{p} \quad(p\gg1),
\end{equation}
with unkown coefficients $c_k$. In \cref{sec:asymp}, we reveal how $c_k$ depend on the geometry of a polygonal domain. 

Second, we look at the impact of the domain symmetry onto the coefficients 
\begin{equation}\label{eq:Ak}
A_k^{(p)} = \frac{1}{\sqrt{|\pa|}} \int_{\pa} v_k^{(p)} d\x,
\end{equation}
where $|\pa|$ is the Lebesgue measure of $\pa$. Note that, integrating \cref{prob:VkSteklov} over $\Omega$ and using the Green's formula, one can also represent $A_k^{(p)}$ as 
\begin{equation}\label{eq:Ak2}
A_k^{(p)} = \frac{p}{\mu_k^{(p)}\sqrt{|\pa|}} \int_{\Omega} V_k^{(p)} d\x.
\end{equation}
These coefficients play an important role in many spectral expansions (see \cite{grebenkov2020paradigm} for details). When $p=0$, one has $\mu_0^{(0)}=0$ and $v_0^{(0)}=1/\sqrt{|\pa|}$ that implies $A_k^{(0)} = \delta_{k,0}$ due to orthogonality of eigenfunctions $v_k^{(0)}$ to $v_0^{(0)}$, where $\delta$ is the Kronecker symbol. Moreover, for any $p\geq0$, one also gets $A_k^{(p)} = \delta_{k,0}$ in the case of a disk due to its rotational symmetry. In \cref{sec:trunc}, we study how the coefficients $A_k^{(p)}$ depend on $p$ and $k$ for various planar shapes, and discuss implications.

Third, we analyze the behavior of the Steklov eigenfunctions $V_k^{(p)}$ away from the boundary. In \cite{hislop2001spectral,polterovich2019nodal,galkowski2019pointwise,daude2021exponential,helffer2022semi}, the localization of $V_k^{(p)}$ near the boundary $\pa$ and exponentially decaying upper bounds were shown for bounded domains with real-analytic boundary $\pa$. In \cref{sec:asymp2}, we inspect the exponential decay of Steklov eigenfunctions away from the boundary for smooth and polygonal domains. In particular, we highlight the role of eigenvalues $\mu_k^{(p)}$ as the decay rates.

\begin{figure}[!ht]
\centering
\subfloat[\label{sfig:disk}]{\includegraphics[width=0.33\linewidth]{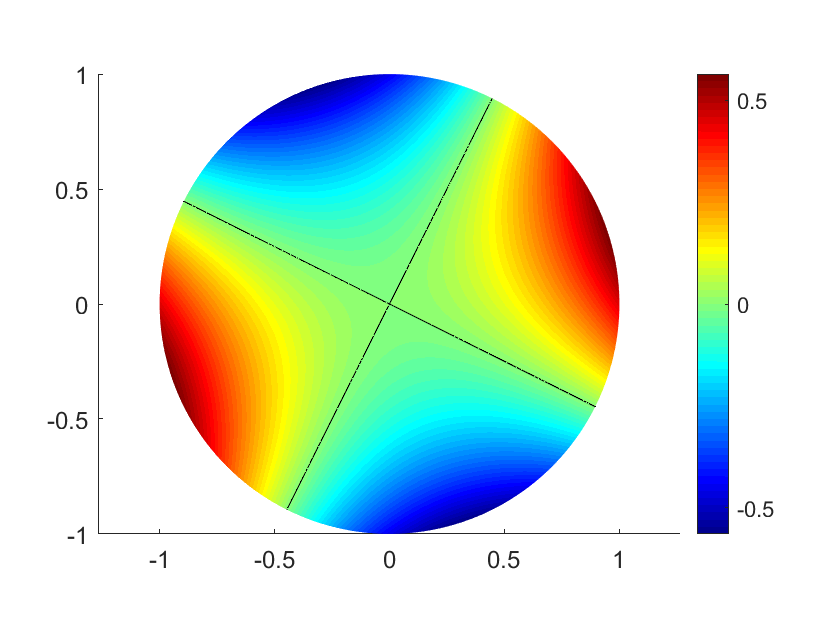}}
\subfloat[\label{sfig:ellip}]{\includegraphics[width=0.33\linewidth]{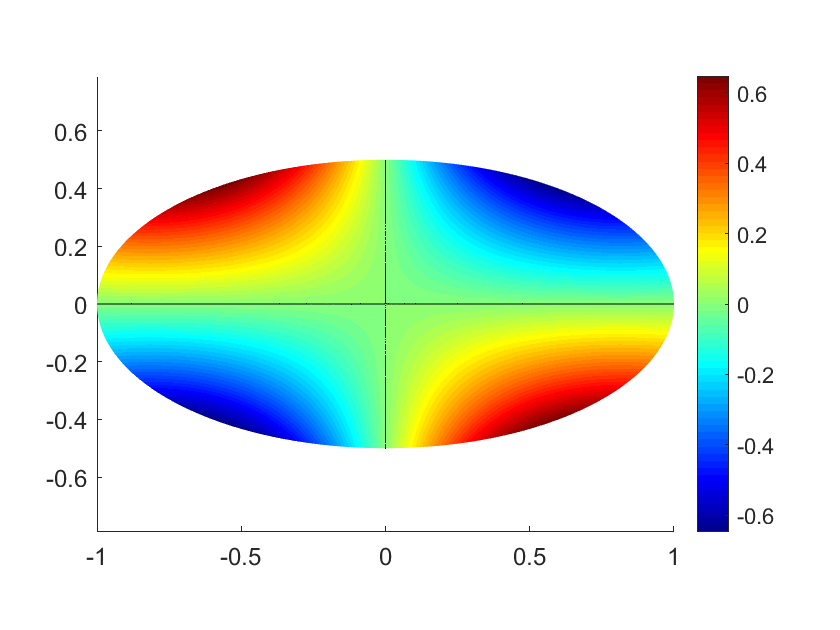}}
\subfloat[\label{sfig:arb}]{\includegraphics[width=0.33\linewidth]{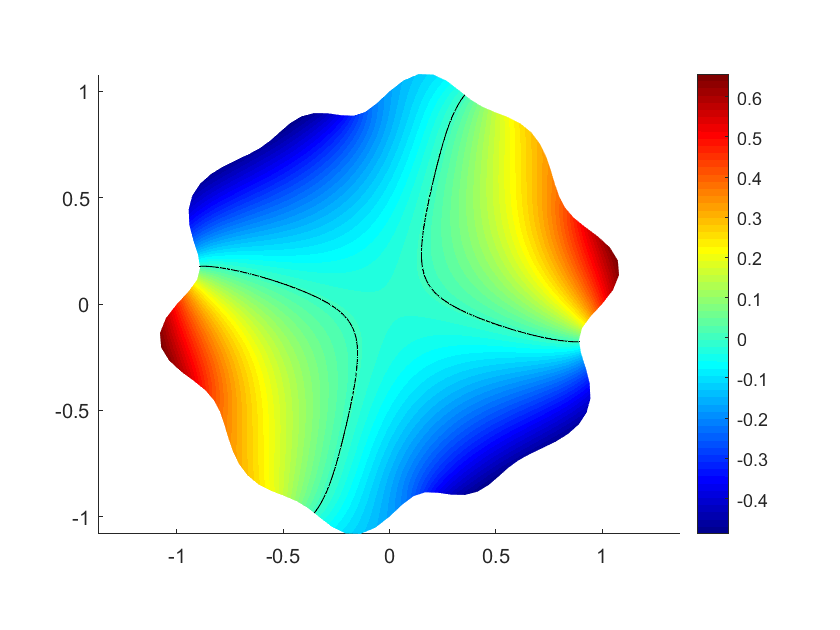}}
\hfil
\subfloat[\label{sfig:tri}]{\includegraphics[width=0.33\linewidth]{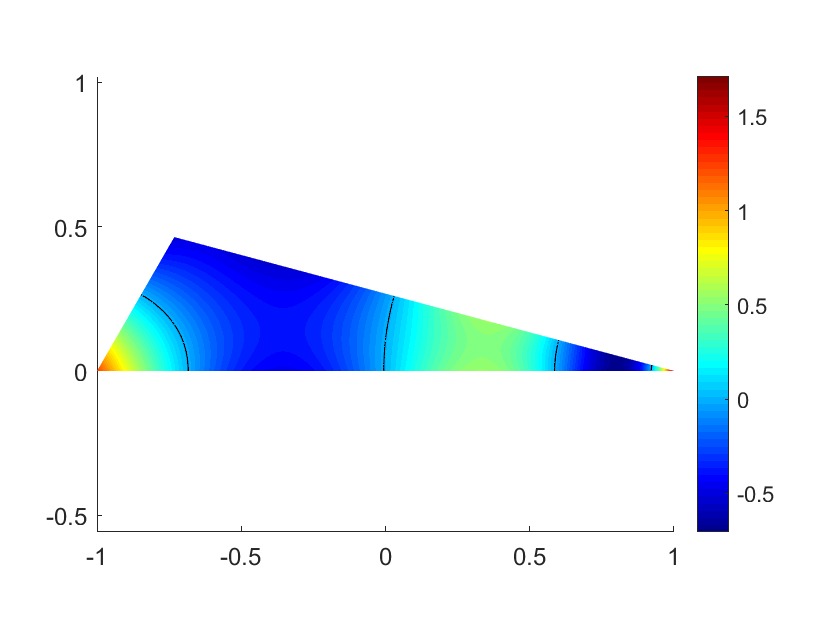}}
\subfloat[\label{sfig:rect}]{\includegraphics[width=0.33\linewidth]{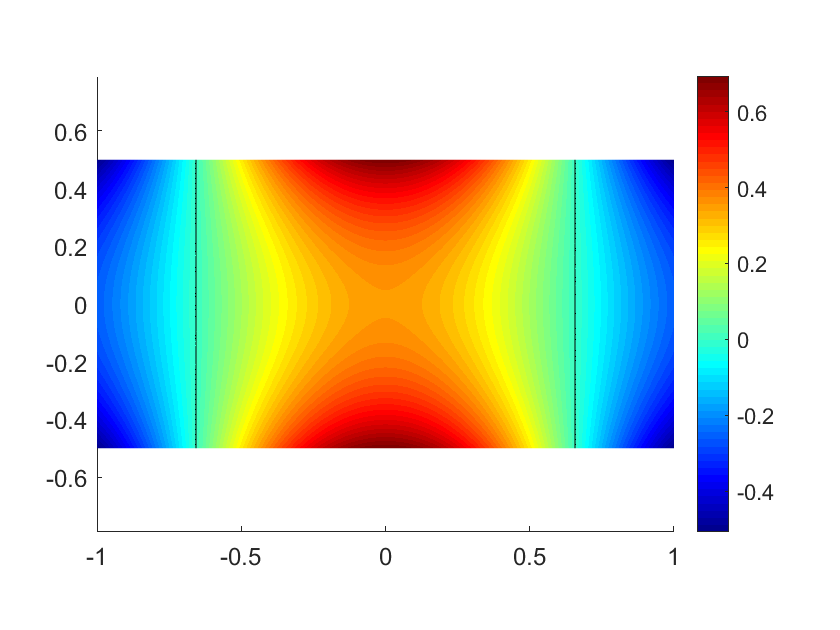}}
\subfloat[\label{sfig:poly8}]{\includegraphics[width=0.33\linewidth]{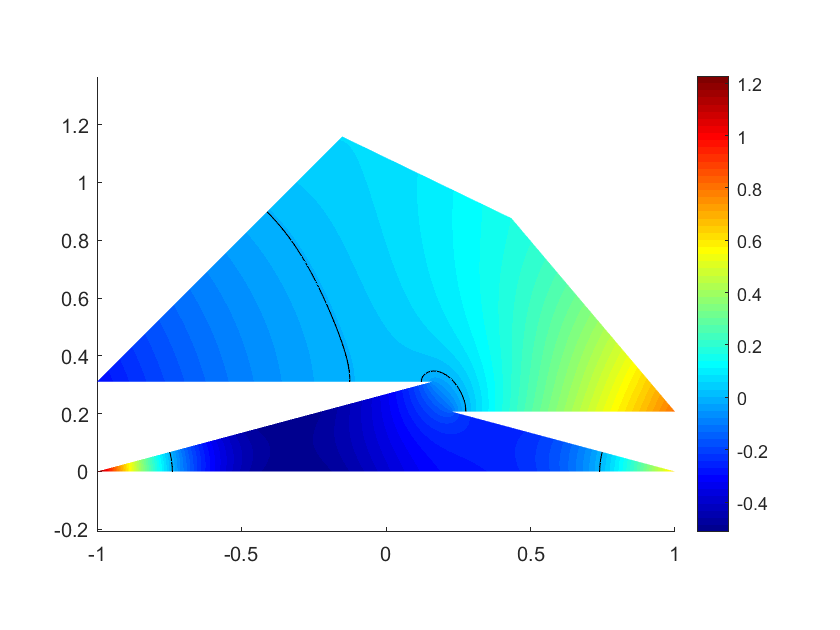}}
\hfil
\subfloat[\label{sfig:frac0}]{\includegraphics[width=0.33\linewidth]{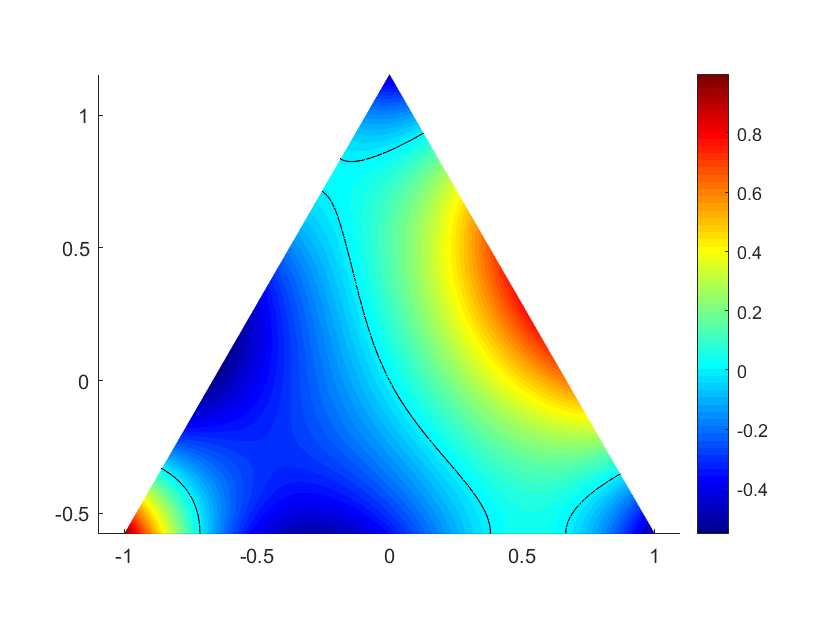}}
\subfloat[]{\includegraphics[width=0.33\linewidth]{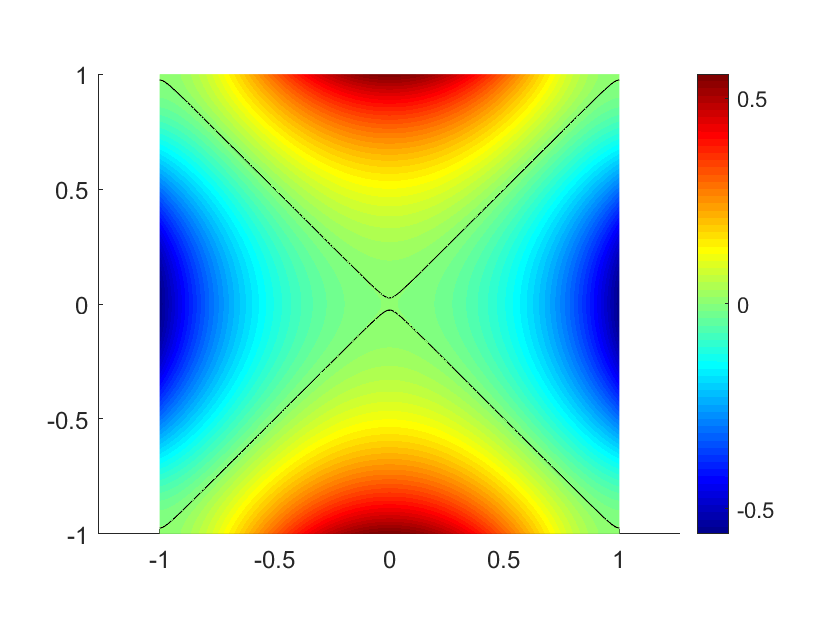}}
\subfloat[\label{sfig:poly5}]{\includegraphics[width=0.33\linewidth]{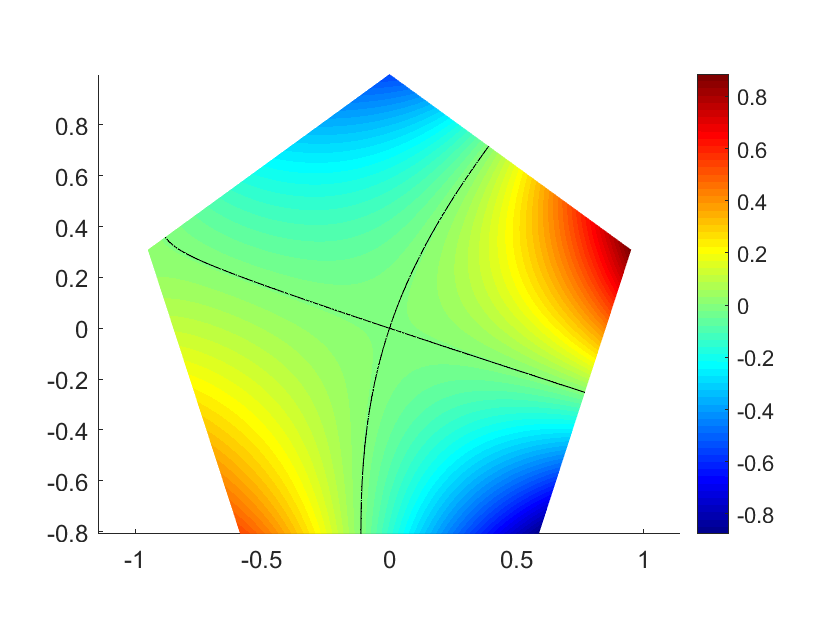}}
\hfil
\subfloat[\label{sfig:frac1}]{\includegraphics[width=0.33\linewidth]{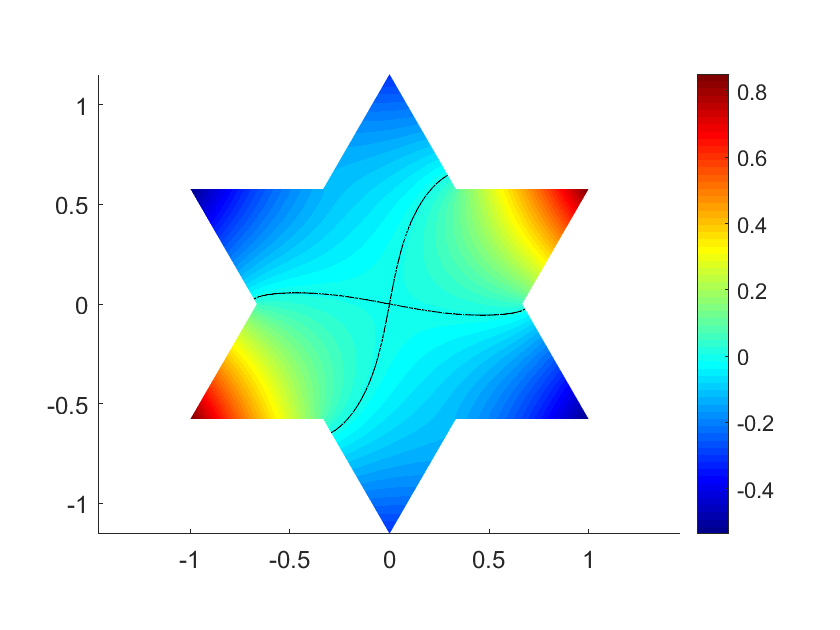}}
\subfloat[\label{sfig:frac2}]{\includegraphics[width=0.33\linewidth]{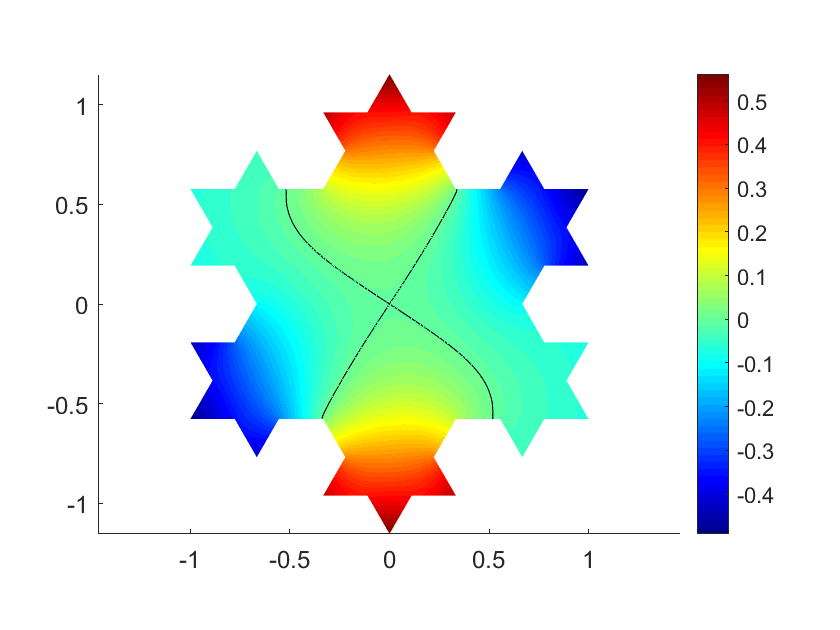}}
\subfloat[\label{sfig:frac3}]{\includegraphics[width=0.33\linewidth]{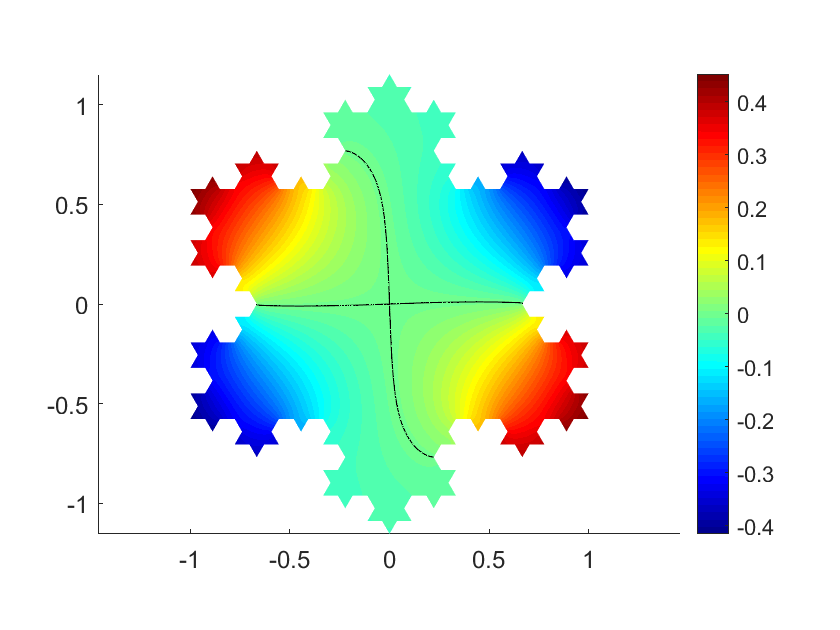}}
\caption{Various planar domains considered in this study and the eigenfunction $V_4^{(1)}$ for each of them. Solid black lines indicate nodal lines.}
\label{fig:Vk}
\end{figure}

To address these questions, we compute $\mu_{k}^{(p)}$, $v_k^{(p)}$ and $V_k^{(p)}$ numerically for a variety of planar domains such as ellipses, triangles, rectangles, regular polygons, Koch snowflakes and randomly generated smooth shapes (see \cref{fig:Vk}).

\section{Numerical methods}\label{sec:nummet}
There is a very limited number of domains whose symmetries allow for separation of variables and thus lead to fully explicit formulas for the eigenvalues and eigenfunctions of $\Mp$ \cite{levitin2023topics,grebenkov2020surface}. In other cases, one needs to employ numerical methods to construct the Dirichlet-to-Neumannn operator. For instance, one can use finite-difference or finite element methods to discretize the problem on a regular lattice or a mesh and to construct a matrix representing the Dirichlet-to-Neumann operator $\Mp$ that needs to be diagonalized to approximate the eigenvalues $\mu_{k}^{(p)}$ and eigenfunctions $v_k^{(p)}$. Flexibility is one of the advantages of these techniques that can deal in the same way with more general second-order elliptic operators. In turn, the need for mesh construction and large sizes of the matrices to be diagonalized are usual drawbacks. Various improvements have been proposed to overcome these limitations: an isoparametric variant of the finite element method for solving
Steklov eigenvalue problems in $\mathbb{R}^d$ for second-order, self-adjoint, elliptic differential operators
\cite{andreev2004isoparametric}, a two-grid discretization scheme \cite{bi2011two,li2011two}, a finite element multi-scale discretization with
an adaptive algorithm based on the shifted inverse iteration \cite{bi2016adaptive}, an iterative multilevel approach
\cite{xie2014type}, a nonconforming finite element methods \cite{li2013nonconforming,yang2009nonconforming}. Another possibility is the reformulation of the Steklov eigenvalue problem in terms of an equivalent boundary integral equation \cite{akhmetgaliyev2017computational, bruno2020domains,chen2020analytical}, or the method of fundamental solutions for solving the Steklov and related spectral problems for the Laplace operator \cite{kupradze1964method,bogosel2016method}. In the planar case, one can also employ conformal mapping to transform the original domain into a simpler domain (e.g., a disk), at the price of dealing with generalized Robin boundary condition \cite{alhejaili2019numerical}. 

For the purpose of our study, the basic finite element method provided a sufficient accuracy and moderate computational cost. Its practical implementation is detailed in \cref{sec:method1}. In addition, we discuss an alternative technique based on the restriction of Green's functions (\cref{sec:method2}).
\subsection{Finite element method}\label{sec:method1}

The first numerical method aims at representing the Dirichlet-to-Neumann operator by a matrix in
two steps. First, we construct the vector $U$ representing the solution $u$ of the modified Helmholtz problem \cref{prob:helm} with a finite element method. 
We discretize $\Omega$ into a triangular mesh, which has $N_i$ nodes inside the domain $\Omega$, and $N_e$ nodes on the boundary $\pa$. The interior nodes are enumerated by $i=1, \ldots, N_i$, while the boundary nodes are enumerated by $i=N_i+1,\ldots , N_i + N_e$. This discretization is equivalent to projection of the weak form of the equation onto a finite-dimensional subspace of dimension $N_p=N_i+N_e$. Let $\{\phi_i\}$, with $i=1,\dots,N_p$, be piecewise polynomial basis functions of a subspace of $\mathcal{H}^1(\Omega)$. In our implementation, each $\phi_i$ is a ``hat" function that is linear on each element and takes the value $0$ at all nodes $x_{j}$ except for $x_i$ at which it is equal to $1$. This property ensures that the functions $\{\phi_i\}_{i=1,\ldots,N_i}$ vanish on the boundary and thus can serve as a basis of $\mathcal{H}^1_0(\Omega)$.

Multiplying the modified Helmholtz equation $(p-\Delta)u=0$ by a test function $\phi_j$ with $j = 1,2,...,N_i$, integrating over $\Omega$, and using the Green's formula, one gets
\begin{equation}\label{eq:weak}
\forall j=1,2,...,N_i, \quad \int_{\Omega} p u \phi_j + \int_{\Omega} \nabla u \nabla \phi_j = 0,
\end{equation}
without boundary terms because $\phi_j|_{\pa}=0$. Next, one can approximate the solution $u$ as a linear combination of basis functions:
\begin{equation}\label{eq:u}
u = \sum_{i=1}^{N_p} U_i \phi_i ,
\end{equation}
with unknown coefficients $U_i$. 
We split the sum as
\begin{equation}\label{eq:split}
u = \sum_{i=1}^{N_i} U_i \phi_i + \underbrace{\sum_{i=N_i+1}^{N_i+N_e} F_i \phi_i}_{=f},
\end{equation}
where the first term represents $u$ inside $\Omega$ and the second one incorporates the Dirichlet boundary condition $u = f$ by setting 
\begin{equation}\label{eq:ue}
\forall i = N_i+1,\ldots,N_i+N_e, \quad U_i = F_i,
\end{equation}
with the coefficients $F_i$ representing $f$. Substitution of this expansion into \cref{eq:weak} yields for all $j=1,2,...,N_i$
\begin{equation}
\sum_{i=1}^{N_i}\left(p \int_{\Omega} \phi_i \phi_j+ \int_{\Omega} \nabla \phi_i \nabla \phi_j\right) U_i =-\sum_{i=N_i+1}^{N_i + N_e}\left(p \int_{\Omega} \phi_i \phi_j+ \int_{\Omega} \nabla \phi_i \nabla \phi_j\right) F_i.
\end{equation}
We denote by $K$ the stiff matrix (of size $N_p\times N_p$) given by $\int_{\Omega} \nabla \phi_i \nabla \phi_j$ and by $M$ the mass matrix (of size $N_p\times N_p$) given by $\int_{\Omega} \phi_i \phi_j$. Since the basis functions $\phi_i$ vanish on all the elements that do not contain the node $x_i$, $K_{ij}$ and $M_{ij}$ are zero except if $x_i$ and $x_j$ are the vertices of the same element and thus the matrices $K$ and $M$ are very sparse.  We get the matrix formulation:
\begin{equation}
\begin{aligned}
(pM+K)^{ii}~U^i &= -(pM+K)^{ie}~F,
\end{aligned}
\end{equation}
where $(pM+K)^{ii}$ is the $N_i \times N_i$ submatrix, $(pM+K)^{ie}$ is the $N_i \times N_e$ submatrix, $F$ is the $N_e \times 1$ vector,  and $U^i$ is the $N_i \times 1$ vector that can thus be found as
\begin{equation}
U^i = -[(pM+K)^{ii}]^{-1} (pM+K)^{ie}~F.
\end{equation}
Combining these coefficients with \cref{eq:ue}, we construct the whole vector $U$ of coefficients $U_i$  as
\begin{equation}\label{eq:U}
U 
= \begin{pmatrix}
-[(pM+K)^{ii}]^{-1} (pM+K)^{ie} \\
I^{ee}
\end{pmatrix}F,
\end{equation}
where $I^{ee}$ is the identity matrix of size $N_e\times N_e$. This concludes the first step.

The second step consists in representing the action of the normal derivative. We restart from the weak formulation of $(p-\Delta)u = 0$
\begin{equation}\label{eq:fvar}
\forall v \in \mathcal{H}^1(\Omega), \quad \int_{\Omega} p u v+\int_{\Omega}  \nabla u \nabla v=\int_{\partial \Omega}  h v,
\end{equation}
with $h=\partial_n u$.
Substituting \cref{eq:u} into \cref{eq:fvar}, the weak formulation reads for any test function $\phi_j$ as

\begin{equation}
\begin{aligned}
&\sum_{i=1}^{N_p}\left(p \int_{\Omega} \phi_i \phi_j+ \int_{\Omega} \nabla \phi_i \nabla \phi_j\right) U_i = \sum_{i=N_i+1}^{N_p}\left( \int_{\partial \Omega} \phi_i \phi_j\right) H_i \quad (j=1,...,N_i),
\end{aligned}
\end{equation}
where $H_i$ are the unknown coefficients representing $h$ on the basis $\{\phi_i\}$.
Denoting by $M_b$ the matrix (of size $N_e \times N_e$) given by $\int_{\partial \Omega} \phi_i \phi_j$, we get the matrix formulation
\begin{equation}
\left( p M+K\right)^{ep} U= M_b H,
\end{equation}
where $\left( p M+K\right)^{ep}$ is the $N_e\times N_p$ submatrix. 
As a consequence, \cref{eq:U} for $U$ implies 
\begin{align}\label{eq:H}
H &= {\bf M}_p F,
\end{align}
where
\begin{equation}
{\bf M}_p = M_b^{-1}\left(p M+K\right)^{ep} \begin{pmatrix}
-[(pM+K)^{ii}]^{-1} (pM+K)^{ie} \\
I^{ee}
\end{pmatrix}.
\end{equation}
According to \cref{eq:H} the matrix ${\bf M}_p$ transforms the Dirichlet boundary condition $u=f$, with a function $f$ represented by the vector $F$, into the Neumann boundary condition $\partial_n u = h$, with $h$ represented by the vector $H$. In other words, this is a matrix representation of the Dirichlet-to-Neumann operator $\mathcal{M}_p$ in terms of basis functions $\{\phi_i\}$. Once the matrix ${\bf M}_p$ is constructed, one can apply standard numerical algorithms to  diagonalize it. The obtained eigenvalues of ${\bf M}_p$ approximate the eigenvalues $\mu_k^{(p)}$. As ${\bf M}_p$ is a finite-size matrix  (of size $N_e \times N_e$), only a finite number of eigenvalues $\mu_k^{(p)}$ can be accurately approximated. In practice, we will limit our analysis to few tens of eigenvalues (say with $k$ from 0 up to 20), for which the method is very accurate, as checked below.
In turn, each eigenvector ${\bf V}_k$ of the matrix ${\bf M}_p$ determines the coefficients of the expansion of $v_k^{(p)}$ on the basis functions:
\begin{equation}
v_k^{(p)}(\x) = \sum\limits_{i=1}^{N_e} {\bf V}_{k,i}~ \phi_{i+N_i}(\x).
\end{equation}
As $\phi_i$ are chosen to be the hat functions, one simply has $v_k^{(p)}(\x_{j+N_i}) = {\bf V}_{k,j}$ at the nodes $\x_{j+N_i}$ of the boundary. As eigenvectors $\bf V_k$ and thus the eigenfunctions $v_k^{(p)}$ are defined up to a multiplicative factor, we explicitly renormalize them to ensure the unit $L^2(\partial\Omega)$ norm, i.e. $\int_{\pa} |v_k^{(p)}|^2 = 1$. 
The Steklov eigenfunctions $V_k^{(p)}$ are obtained by replacing $F$ by $v_k^{(p)}$ in \cref{eq:U}. 

The flexibility of the finite element method allows for various extensions. For instance, one can solve the mixed Steklov problem when the boundary $\pa$ is composed of two disjoint parts: $\pa=\pa_1 \cup \pa_2$. The Dirichlet-to-Neumann operator then acts as 
\begin{equation}
\begin{aligned}
\Mp: H^{1 / 2}(\pa_1) &\rightarrow H^{-1 / 2}(\pa_1) \\
f &\mapsto \left.\left(\partial_n u\right)\right|_{\pa_1},
\end{aligned}
\end{equation}
where 

\begin{minipage}{.5\linewidth}
\begin{equation*}
\left\{ \begin{array}{r l l}
(p-\Delta) u(\x)& = 0 & \quad (\x \in \Omega), \\
u(\x)& = f(\x) &\quad  (\x \in \pa_1), \\
u(\x)& = 0 &\quad  (\x \in \pa_2),
\end{array}  \right.
\end{equation*}
\end{minipage}%
or
\begin{minipage}{.5\linewidth}
\begin{equation*}
\left\{ \begin{array}{r l l}
(p-\Delta) u(\x)& = 0 & \quad (\x \in \Omega), \\
u(\x)& = f(\x) &\quad  (\x \in \pa_1), \\
\partial_n u(\x)& = 0 &\quad  (\x \in \pa_2).
\end{array}  \right.
\end{equation*}
\end{minipage} \\

\noindent In other words, one imposes an additional Dirichlet or Neumann boundary condition on $\pa_2$ and defines the operator $\Mp$ to act on functions on $\pa_1$.
In the context of diffusion-controlled reactions, mixed boundary conditions allow one to describe various processes in which the diffusing particle can leave
the confining domain through an escape region or be destroyed on it, before reaching the target region \cite{grebenkov2023encounter}.  The implementation of the Dirichlet boundary condition to our numerical method is straightforward: one just needs to truncate the second term in \cref{eq:split} to the indices corresponding to the nodes on $\pa_1$. In this way, one incorporates the Dirichlet boundary condition $u=f$ on $\pa_1$, while letting the homogeneous Dirichlet boundary condition $u=0$ on $\pa_2$. The inclusion of the Neumann boundary condition is as well simple: one extends the first term in \cref{eq:split} to the indices corresponding to the nodes on $\pa_2$, i.e. the first term then represents $u$ inside $\Omega$ and $u$ on the reflecting boundary $\pa_2$, while the second term (reduced to the indices corresponding to the nodes on $\pa_1$) encodes the Dirichlet boundary condition $u=f$ on $\pa_1$. An extension of this method to three-dimensional domains is also straightforward.

\subsection{Green's function method} \label{sec:method2}
The second method was inspired by the spectral decompositions of the Green's functions \cite{grebenkov2020paradigm}. Let us introduce the Green's function $\tilde{G}_q(\x,p|\x_0)$ satisfying
\begin{equation}\label{prob:G}
\left\{ \begin{array}{r l l}
(p-\Delta) \tilde{G}_q(\x,p|\x_0)  & = & \delta(\x-\x_0) \quad (\x \in \Omega), \\
-\partial_n \tilde{G}_q(\x,p|\x_0)& = &q \tilde{G}_q(\x,p|\x_0) \quad (\x \in \partial \Omega),
\end{array}  \right.
\end{equation}
with a constant $0 \leq q \leq \infty$, and $\delta(\x-\x_0)$ being the Dirac distribution.
Its expansion on the Steklov eigenbasis reads \cite{grebenkov2020paradigm}:
\begin{equation}\label{eq:decspecDtN1}
\tilde{G}_q(\x,p|\x_0) = \tilde{G}_\infty(\x,p|\x_0) + \sum_{k=0}^{\infty} \frac{[V_k^{(p)}(\x_0)]^* V_k^{(p)}(\x) }{\mu_k^{(p)} + q},
\end{equation}
which is also possible to write as
\begin{equation}\label{eq:decspecDtN2}
\tilde{G}_q(\x,p|\x_0) = \tilde{G}_0(\x,p|\x_0) - \sum_{k=0}^{\infty} \frac{[V_k^{(p)}(\x_0)]^* V_k^{(p)}(\x) }{(\mu_k^{(p)}/q + 1)\mu_{k}^{(p)}},
\end{equation}
where asterisk denotes the complex conjugate. Setting $q = 0$ and restricting the points $\x$ and $\x_0$ onto the
boundary $\pa$ in \cref{eq:decspecDtN1}, one gets the integral kernel of the inverse of the
Dirichlet-to-Neumann operator:
\begin{equation}\label{eq:G0}
\tilde{G}_0(\x,p|\x_0) =  \sum_{k=0}^{\infty} \frac{[v_k^{(p)}(\x_0)]^* v_k^{(p)}(\x) }{\mu_k^{(p)}} \quad (\x, \x_0 \in \pa).
\end{equation}
However, the Green's function $\tilde{G}_0(\x,p|\x_0)$ exhibits a
singularity at $\x = \x_0$ that would require additional regularization  and
may enhance numerical errors.  It is therefore convenient to remove this
singularity by using \cref{eq:decspecDtN2}.
Restricting both $\x$ and $\x_0$ to the boundary, we consider the integral kernel:
\begin{equation}  
\nonumber \tilde{g}_q(\x,p|\x_0) = \frac{\tilde{G}_0(\x,p|\x_0) - \tilde{G}_q(\x,p|\x_0)}{q} = \sum_{k=0}^{\infty} \frac{[v_k^{(p)}(\x_0)]^* v_k^{(p)}(\x) }{\mu_k^{(p)}(\mu_k^{(p)} + q)} \quad (\x, \x_0 \in \Omega).
\end{equation}
One gets then
\begin{equation}
\int_{\pa}  \tilde{g}_q(\x,p|\x_0) v_k^{(p)}(\x) d\x = \eta_k^{(p)} v_k^{(p)}(\x_0),
\end{equation}
where the eigenvalues $\eta_k ^{(p)}$ of this integral operator are related to
$\mu_k^{(p)}$ as
\begin{equation}
\mu_k^{(p)} = \sqrt{1/\eta_k^{(p)} + q^2/4} - q/2 .
\end{equation}
The discretization of the integral at boundary points $\x_i$ yields a
system of linear equations
\begin{equation}
\sum\limits_j \delta\x_j \, \tilde{g}_q(\x_j,p|\x_i)\,  v_k^{(p)}(\x_j) = \eta_k^{(p)} v_k^{(p)}(\x_i),
\end{equation}
where $\delta\x_j$ are the areas of the surface elements around $\x_j$.  In
other words, one needs to diagonalize the matrix $G_{ij} = \delta\x_j
\, \tilde{g}_q(\x_j,p|\x_i)$. 

Each Green's function can be found either directly (e.g., by a finite element method), or from
its (truncated) spectral expansion
\begin{equation}\label{eq:34}
\tilde{G}_q(\x,p|\x_0) = \sum_{k=1}^{\infty} \frac{[u_k^{(q)}(\x_0)]^* u_k^{(q)}(\x)}{p+\lambda_k^{(q)}},
\end{equation}
where $\lambda_k^{(q)}$ and $u_k^{(q)}$ are the eigenvalues and the $L^2(\Omega)$-normalized eigenfunctions of the Laplace operator:
\begin{equation}\label{prob:lap}
\left\{ \begin{array}{r l l}
-\Delta u_k^{(q)}  = \lambda_k^{(q)} u_k^{(q)} \quad (\x \in \Omega), \\
-\partial_n u_k^{(q)} = q u_k^{(q)} \quad (\x \in \partial \Omega).
\end{array}  \right.
\end{equation}
We stress that the use of the spectral expansion \cref{eq:34} is not optimal from the numerical point of view. In fact, the construction of Laplacian eigenfunctions and eigenvalues is a time-consuming procedure, while the truncation of the infinite series in \cref{eq:34} to a finite number of terms can be the major source of numerical errors. We used this expansion as a straightforward way to access the Green's function for validation purposes, but more efficient numerical tools can be designed for this task. In turn, the advantage of this method is that the Laplacian
eigenfunctions need to be computed only once for a given domain.

Once $\mu_k^{(p)}$ and $v_k^{(p)}$ are found, one can also compute the extension $V_k^{(p)}(\x)$ by using the Green's
function and the Neumann boundary condition $\partial_n V_k^{(p)}
= \mu_k^{(p)} v_k^{(p)}$, so that
\begin{equation}  \label{eq:Vk_G0}
V_k^{(p)}(\x_0) = \int_{\pa}\tilde{G}_0(\x,p|\x_0) \,  \mu_k^{(p)} v_k^{(p)}(\x) d\x.
\end{equation}
Note that this relation is not applicable for $p = 0$ and $k=0$, for
which $\mu_0^{(0)} = 0$ and $V_0^{(0)}(\x_0) = 1/\sqrt{|\pa|}$ are already known.

\subsection{Numerical validation}

Both methods were implemented in \textit{Matlab}. In particular, we relied on the \textit{Matlab PDEtool} to generate triangular meshes, the matrices $K$, $M$, and $M_b$ in the first method, and the Laplacian eigenfunctions in the second method. We  also used build-in functions \textit{eig} and \textit{eigs} for matrix diagonalizations.

We validate both numerical methods by comparing their results with the explicit formulas known for the disk of radius $R$ \cite{grebenkov2020surface}: 
\begin{subequations}\label{eq:eigelem}
\begin{align}
\mu_{0}^{(p)} &= \sqrt{p} \frac{I_1(R \sqrt{p})}{I_0(R\sqrt{p})}, \quad v_0^{(p)}(\theta) = \frac{1}{\sqrt{2\pi R}}, \quad V_0^{(p)}(r,\theta) = \frac{I_0(r\sqrt{p})}{I_0(R\sqrt{p})} v_0^{(p)}(\theta), \\
\mu_{2k}^{(p)} &= \sqrt{p} \frac{I_{2k}'(R \sqrt{p})}{I_{2k} (R\sqrt{p})}, \quad v_{2k}^{(p)}(\theta) = \frac{\cos(k\theta)}{\sqrt{\pi R}}, \quad V_{2k}^{(p)}(r,\theta) = \frac{I_{2k}(r\sqrt{p})}{I_{2k}(R\sqrt{p})}v_{2k}^{(p)}(\theta), \\
\mu_{2k+1}^{(p)} &= \sqrt{p} \frac{I_{2k+1}'(R \sqrt{p})}{I_{2k+1} (R\sqrt{p})}, \quad v_{2k+1}^{(p)}(\theta) = \frac{\sin(k\theta)}{\sqrt{\pi R}}, \quad V_{2k+1}^{(p)}(r,\theta) = \frac{I_{2k+1}(r\sqrt{p})}{I_{2k+1}(R\sqrt{p})}v_{2k+1}^{(p)}(\theta),
\end{align}
\end{subequations}
where $I_k(z)$ is the modified Bessel function of the first kind, prime denotes the derivative with respect to the argument, and we use polar coordinates $(r,\theta)$.
Note that all the eigenvalues are twice degenerate, except for the first one $\mu_0^{(p)}$ which is simple.
For Method 1 the mesh is composed of 384138 triangles, and the maximal mesh size is 0.0042. For Method 2 the mesh is composed of 16256 triangles and the series in \cref{eq:34} was truncated to 131 eigenfunctions of the Laplace operator. \Cref{tab:diskeigv} summarizes the first $11$ eigenvalues of $\Mp$ for $R=1$ and $p=1$, while \cref{sfig:disk} shows the corresponding eigenfunction $V_4^{(1)}$. One sees that the eigenvalues in the third column, which were numerically obtained by Method 1, are in excellent agreement with the exact ones given by \cref{eq:eigelem}. While the numerical eigenvalues computed by Method 2 are less accurate, they were computed much faster on a mesh with a smaller number of triangles. We conclude that two methods provide complementary numerical tools to access the spectral properties of the Dirichlet-to-Neumann operator. As a systematic comparison of two methods is beyond the scope of the paper, we use Method 1 in the following computations. \Cref{tab:diskeigf} presents the root mean squared errors (RMSE) between the exact eigenfunctions from \cref{eq:eigelem} and the numerical ones. For each index $k$ we compute the RMSE from the formula: $\sqrt{N^{-1} \sum_{i=1}^{N} (f_i-f_i^*)^2}$, where $N$ is the number of points on the boundary, and $f$ and $f^*$ represent the analytical and numerical values. Expectedly, this error increases with the index $k$ but remains negligible for the considered range of indices. In Appendix \ref{secap:num}, we provide additional verifications by presenting the explicitly known formulas for a rectangle and comparing them with our numerical results.

In the following, we set the maximal mesh size to $0.005$.	

%

\begin{table}[!ht]
\begin{minipage}{.45\linewidth}
\centering
\medskip
\begin{tabular}[!hb]{|c|c|c|c|}
\hline \textbf{Index k}  & \textbf{Exact} & \textbf {Method 1}  & \textbf {Method 2} \\
\hline 
0 & 0.4464 & 0.4464 & 0.4464 \\
1 & 1.2402 & 1.2402 & 1.2402 \\
2 & 1.2402 & 1.2402 & 1.2402 \\
3 & 2.1633 & 2.1633 & 2.1640 \\
4 & 2.1633 & 2.1633 & 2.1640 \\
5 & 3.1235 & 3.1235 & 3.1257 \\
6 & 3.1235 & 3.1235 & 3.1259 \\
7 & 4.0992 & 4.0993 & 4.1059 \\
8 & 4.0992 & 4.0993 & 4.1061 \\
9 & 5.0828 & 5.0832 & 5.0955 \\
10 & 5.0828 & 5.0832 & 5.0959 \\
\hline
\end{tabular}
\caption{List of the first 11 eigenvalues $\mu_k^{(p)}$ for the unit disk, with $p=1$.}
\label{tab:diskeigv}
\end{minipage}\hfill
\begin{minipage}{.45\linewidth}
\centering
\begin{tabular}[!hb]{|c|c|c|}
\hline \textbf{Index k}  & \textbf{RMSE 1} & \textbf{RMSE 2}  \\
\hline 
0 & 0.0000 & 0.0000 \\
1 & 0.0013 & 0.0004\\
2 & 0.0013 & 0.0004\\
3 & 0.0021 & 0.0006 \\
4 & 0.0021 & 0.0055 \\
5 & 0.0044 & 0.0054 \\
6 & 0.0044 & 0.0054 \\
7 & 0.0024 & 0.0173\\
8 & 0.0024 & 0.0173\\
9 & 0.0031 & 0.0120\\
10 & 0.0031 & 0.0120 \\
\hline
\end{tabular}
\caption{List of root mean squared errors for the first $11$ eigenfunctions $v_k^{(p)}$ of the unit disk, with $p=1$.}
\label{tab:diskeigf}
\end{minipage}
\end{table}

\section{Asymptotic behavior of eigenvalues}\label{sec:asymp}
In this section, we investigate the asymptotic behavior of the eigenvalues $\mu_k^{(p)}$ in the limit $p\to \infty$. In particular, we focus on the dependence of $\mu_k^{(p)}$ on the angles of polygonal domains.

\subsection{Ellipses} 
\begin{figure}[!ht]
\centering
\subfloat[Aspect ratio $2$]{\includegraphics[width=0.455\linewidth]{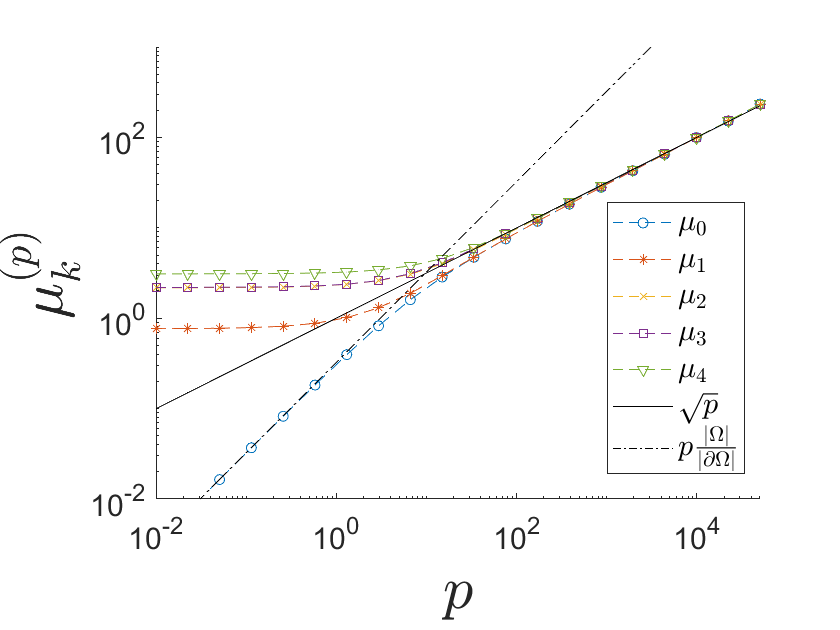}}
\hfil
\subfloat[Aspect ratio $10$\label{sfig:mu-vs-p_ellipN=10}]{\includegraphics[width=0.455\linewidth]{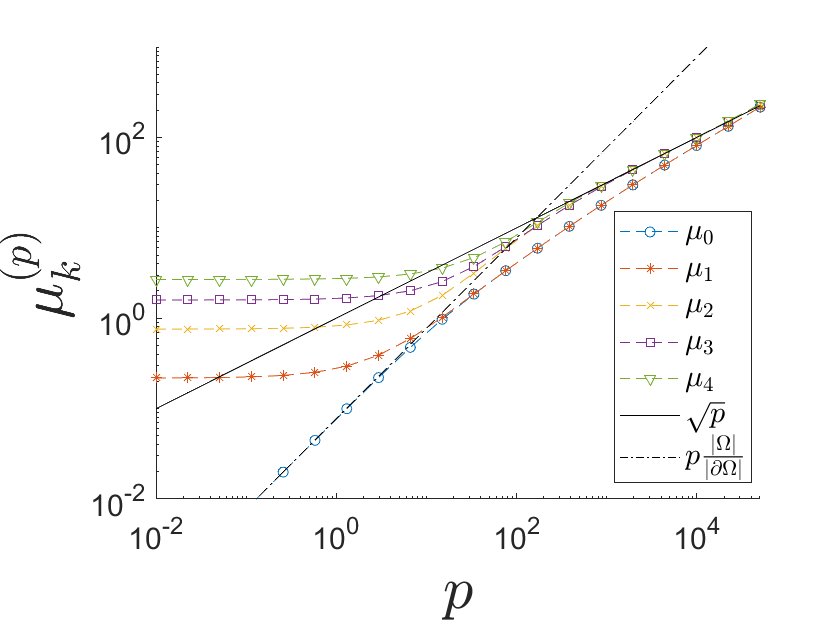}}
\caption{Dependence of the first eigenvalues $\mu_k^{(p)}$ (in symbol) on $p$ for an ellipse with semiaxes 1 and 0.5 (a) and an ellipse with semiaxes 1 and 0.1 (b). Solid black line presents the asymptotic relation (\ref{eq:asymp-mu_circle}) and dotted black line indicates the relation (\ref{eq:asymp-mu2}).}
\label{fig:mu-vs-p_ellip}
\end{figure}

We start the study by considering smooth anisotropic domains such as an ellipse with semiaxes $a$ and $b$ (\cref{sfig:ellip}):
\begin{equation}
\Omega= \{ (x,y) \in \mathbb{R}^2 : (x/a)^2 + (y/b)^2 <1\}.
\end{equation}
In this configuration, we numerically check the expected asymptotic relation \cref{eq:asymp-mu_circle} at large $p$ as well as the opposite limit $p\to0$ \cite{lacey1998multidimensional,levitin2008principal}:
\begin{equation}\label{eq:asymp-mu2}
\mu_0^{(p)} \simeq p \frac{|\Omega|}{|\pa|} \quad(p\ll 1).
\end{equation}
\Cref{fig:mu-vs-p_ellip} illustrates these relations for two ellipses of aspect ratio $2$ and $10$.
For the second (more elongated) ellipse (\cref{sfig:mu-vs-p_ellipN=10}), the first two eigenvalues reach the asymptotic behavior $\sqrt{p}$ slower than the others, showing an extended transient regime at intermediately large $p$. The existence of this transient regime suggests that the dependence of eigenvalues on $p$ may reveal additional geometrical features of the domains, such as its anisotropy.

\subsection{Rectangles and regular polygons}

\begin{figure}[!ht]
\centering
\subfloat[\label{fig:coefa-vs-poly}]{\includegraphics[width=0.455\linewidth]{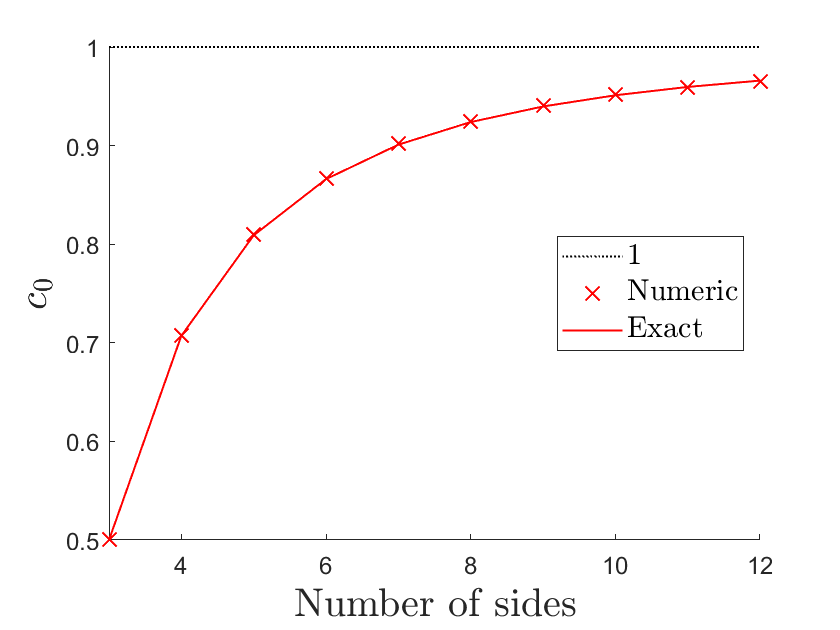}}
\hfil
\subfloat[\label{fig:coefa-vs-rect}]{\includegraphics[width=0.455\linewidth]{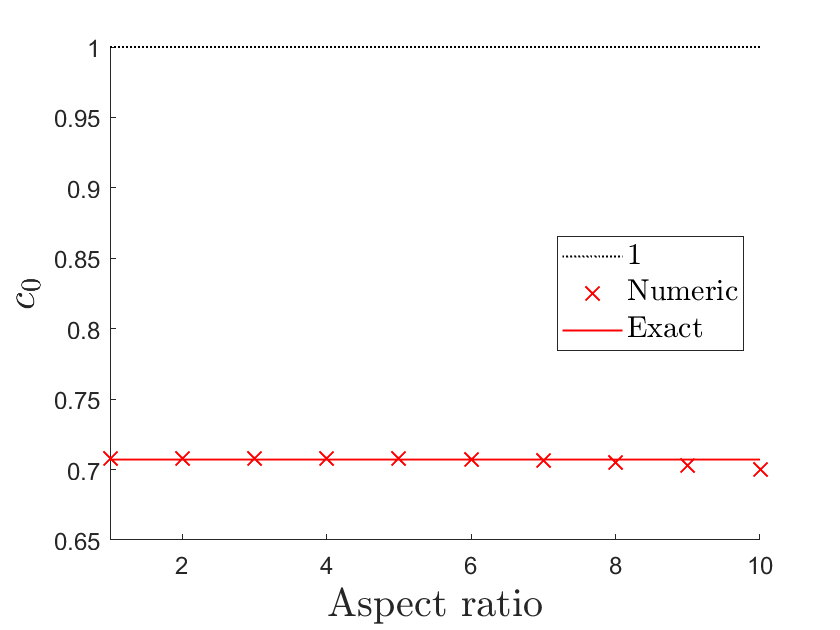}}
\caption{The coefficient $c_0$ (in symbol) as a function of (a) the number $N$ of sides of the regular polygon; (b) the aspect ratio of the rectangle with one side of length 1 and the other varying from 1 to 10. Solid line presents \cref{eq:conj1}, while dashed horizontal line indicates the value 1 of this coefficient for domains with smooth boundaries.}
\end{figure}

Then, we inspect the role of angles of a polygonal domain on the asymptotic behavior of the eigenvalues. For this purpose, we consider rectangular and regular polygonal domains. For a square, as $p \to \infty$, we observe the asymptotic behavior \cref{eq:asymp-ck}, with $c_k \approx 0.51$ for $k=\{0, 1,2,3\}$, and $c_k = 1$ for other $k$. The prefactor $c_k$ was estimated from the ratio $\mu_k^{(p)}/\sqrt{p}$ at $p=10^3$. In other words, the first four eigenvalues exhibit the asymptotic behavior \cref{eq:asymp-ck} and deviate from  \cref{eq:asymp-mu_circle}, which was demonstrated for bounded domains with smooth boundaries. Note that $\mu_1^{(p)}$ and $\mu_2^{(p)}$ are identical, i.e., they are degenerate eigenvalues of multiplicity 2. 
%
%
In the case of a regular polygon with $N$ vertices, the first $N$ eigenvalues exhibit the asymptotic behavior \cref{eq:asymp-ck}, with the same coefficients $c_0 =... =c_{N-1} < 1$. As the number of sides of the polygon increases, the domain is getting closer to a disk so that the coefficients $c_0=...=c_{N-1}$ approach $1$ (see \cref{fig:coefa-vs-poly}). We conjecture that these coefficients only depend on the angle $\alpha=\pi(1-2/N)$ of the regular polygon and take the value 
\begin{equation}\label{eq:conj1}
c_0=...=c_{N-1} = \sin(\alpha/2).
\end{equation} 
This conjecture is numerically confirmed on \cref{fig:coefa-vs-poly}.

To underline the role of the angles, we plot the coefficient $c_0$ as a function of the aspect ratio of a rectangle (\cref{fig:coefa-vs-rect}). It appears that the anisotropy of the domain does not affect the coefficient $c_0$, whose minor variations can be attributed to weak inaccuracies of the numerical method for elongated domains.

\subsection{Prefractals}

\begin{figure}[!ht]
\centering
\subfloat[Generation 0 (equilateral triangle) \label{fig:mu-vs-p_fractal0}]{\includegraphics[width=0.45\linewidth]{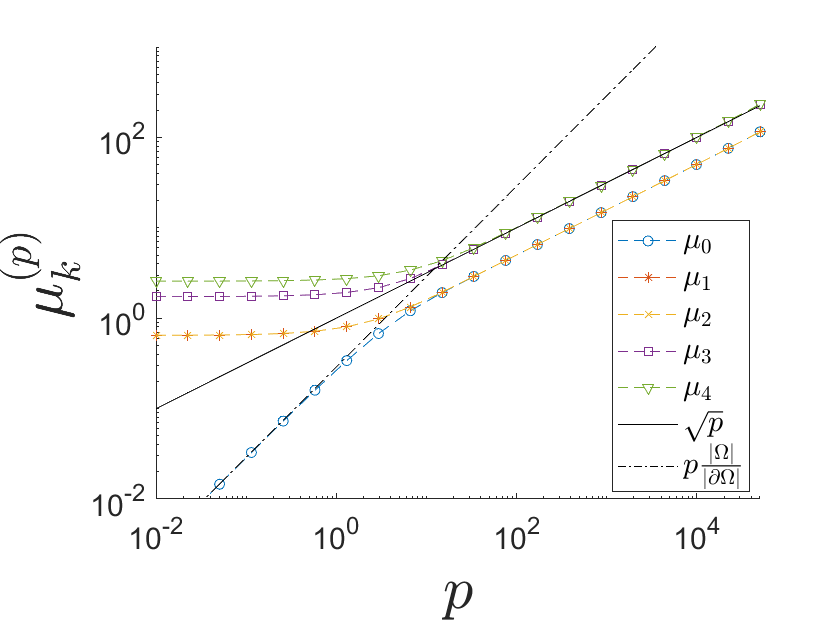}}
\hfil
\subfloat[Generation 1]{\includegraphics[width=0.45\linewidth]{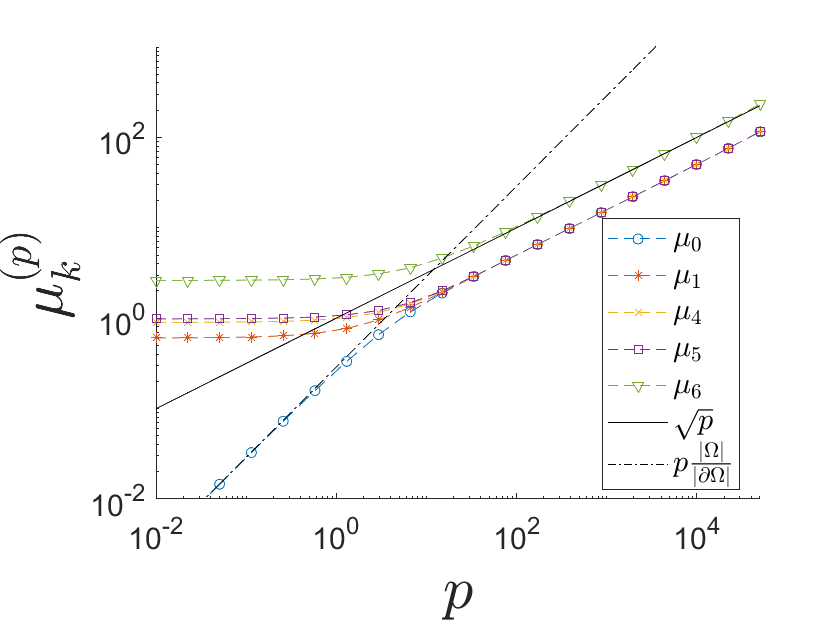}}
\hfil
\subfloat[Generation 2]{\includegraphics[width=0.45\linewidth]{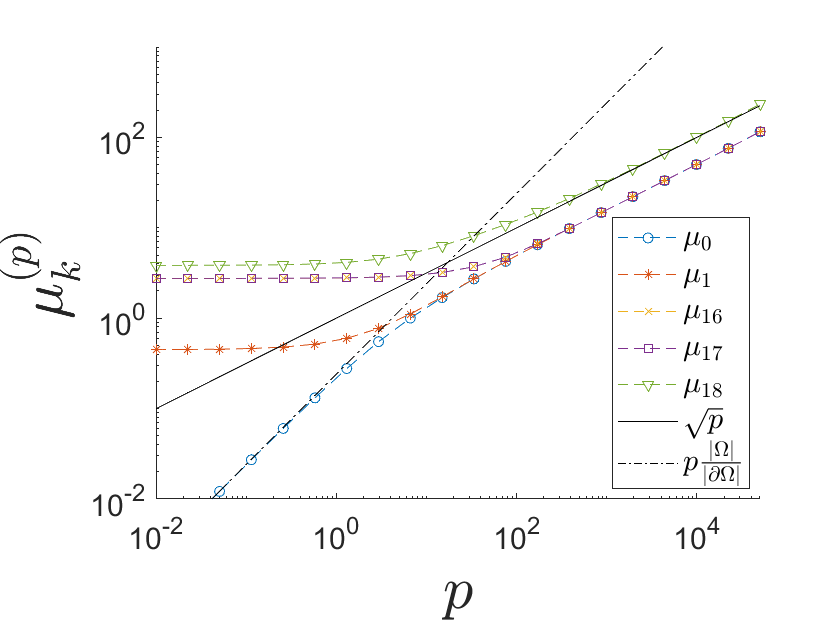}}
\hfil
\subfloat[Generation 3 \label{fig:mu-vs-p_fractal3}]{\includegraphics[width=0.45\linewidth]{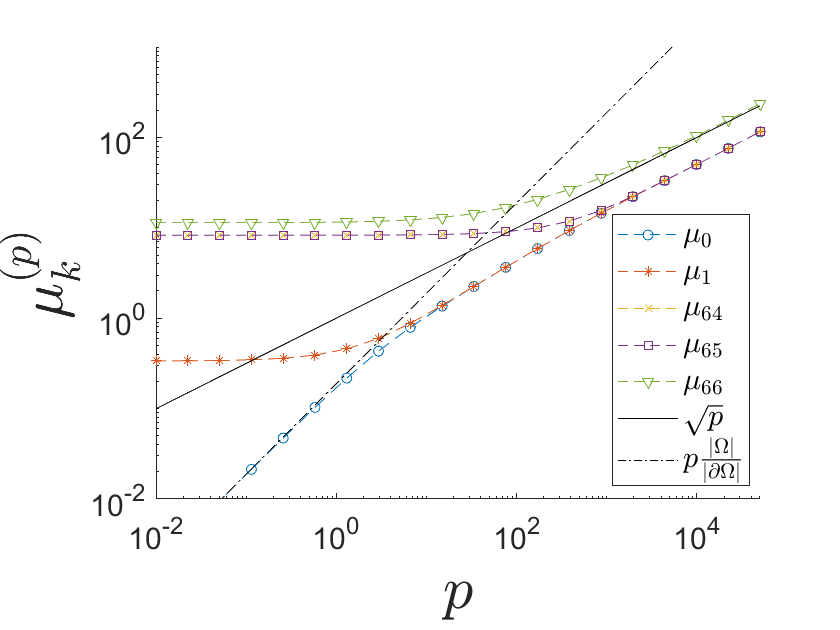}}
\caption{Dependence of some eigenvalues $\mu_k^{(p)}$ (in symbol) on $p$ for the first four generations of the Koch snowflake. Solid black line presents the asymptotic relation (\ref{eq:asymp-mu_circle}) and  dotted black line indicates the relation (\ref{eq:asymp-mu2}).}
\label{fig:mu-vs-p_fractal}
\end{figure}

To investigate the effect of roughness of the boundary, we consider a family of prefractal Koch snowflakes. These domains are constructed iteratively, starting from the equilateral triangle of sidelength 2 (generation 0, \cref{sfig:frac0}) and adding finer geometric features at each iteration (see Figs. \ref{sfig:frac1}, \ref{sfig:frac2}, \ref{sfig:frac3} for generations 1, 2 and 3 respectively).
\Cref{fig:mu-vs-p_fractal0} shows the dependence of the eigenvalues on $p$ for an equilateral triangle of sidelength $2$. As $p \to \infty$, we observe the asymptotic behavior \cref{eq:asymp-ck}, with $c_k \approx \sin(\pi/6)=0.5$ for $k=\{0, 1,2\}$, and $c_k = 1$ for other $k$. This is in agreement with the conjectured expression \cref{eq:conj1} of $c_k$ for regular polygons. For next three generations (see \cref{fig:mu-vs-p_fractal}), it appears that the first $6$, $18$ and $66$ eigenvalues respectively have $c_k \approx \sin(\pi/6)=0.5$ for $k$ from 0 to $K_g-1$, and $c_k = 1$ for other $k$,  where $K_g=3(1+4^0+4^1+...+4^{g-1})=2+4^g$, is the number of angles $\pi/3$ in the generation $g$. We conclude that the value of the coefficient $c_0$ is independent of the fractal generation $g$, which only influences the number of $c_k$ that are smaller than 1. Note also that for the considered prefactal domains, the asymptotic behavior \cref{eq:asymp-mu2} still holds, which suggests that the boundary $\pa$ is not rough enough to break this asymptotic relation. 

\subsection{Generic triangle}
We investigate the asymptotic behavior of eigenvalues for a generic triangle, which is constructed by setting the length of one side to be 2 and two angles to be $\pi/12$ and $\pi/3$ (\cref{sfig:tri}). To ensure the quality of the numerical results, we refine the mesh and set the maximal mesh size to 0.003.
%
%
%
The last column of \Cref{tab:ck2} presents the obtained coefficients $c_k$. 
\begin{table}[!ht]
\centering
\renewcommand{\arraystretch}{2.5}
$\begin{array}{|c|ccc|cc|}
\hline \boldsymbol{k}  & \boldsymbol{\alpha_0} & \bf \boldsymbol{\alpha_1} & \boldsymbol{\alpha_2} & \boldsymbol{c_k} \textbf{(conjecture)} & \boldsymbol{c_k} \textbf{(numeric)}  \\
\hline 
0 & \cfrac{\pi}{12} & \cfrac{\pi}{3} & \cfrac{7\pi}{12} &  \sin\left(\cfrac{1}{2}\cdot\cfrac{\pi}{12}\right)=0.1305 &  0.1305 \\
1 & \cfrac{3\pi}{12} & \cfrac{\pi}{3} & \cfrac{7\pi}{12} &  \sin\left(\cfrac{1}{2}\cdot\cfrac{3\pi}{12}\right)=0.3827 &  0.3833 \\
2 & \cfrac{5\pi}{12} & \cfrac{\pi}{3} & \cfrac{7\pi}{12} &  \sin\left(\cfrac{1}{2}\cdot\cfrac{\pi}{3}\right)=0.5 &  0.45999 \\
3 & \cfrac{5\pi}{12} & \textcolor{gray}{\pi} & \cfrac{7\pi}{12} &  \sin\left(\cfrac{1}{2}\cdot\cfrac{5\pi}{12}\right)=0.6088 &  0.6104 \\
4 & \cfrac{7\pi}{12} & \textcolor{gray}{\pi} & \cfrac{7\pi}{12} &  \sin\left(\cfrac{1}{2}\cdot\cfrac{7\pi}{12}\right)=0.7934 &  0.7935 \\
5 & \cfrac{7\pi}{12} & \textcolor{gray}{\pi} & \textcolor{gray}{\cfrac{7\pi}{4}} &  \sin\left(\cfrac{1}{2}\cdot\cfrac{7\pi}{12}\right)=0.7934 &  0.7957 \\
6 & \cfrac{9\pi}{12} & \textcolor{gray}{\pi} & \textcolor{gray}{\cfrac{7\pi}{4}} &  \sin\left(\cfrac{1}{2}\cdot\cfrac{9\pi}{12}\right)=0.9239 &  0.9260 \\
7 & \cfrac{11\pi}{12} & \textcolor{gray}{\pi} & \textcolor{gray}{\cfrac{7\pi}{4}} &  \sin\left(\cfrac{1}{2}\cdot\cfrac{11\pi}{12}\right)=0.9914 &  0.9925 \\
8 & \textcolor{gray}{\cfrac{13\pi}{12}} & \textcolor{gray}{\pi} & \textcolor{gray}{\cfrac{7\pi}{4}} &  \sin\left(\pi\right)=1 &  1.0020 \\
\hline
\end{array}$
\caption{First 9 coefficients $c_k$ for a triangle with angles $\pi/12$, $\pi/3$ and $7\pi/12$. The last column presents the values $c_k$ estimated from the ratio $\mu_k^{(p)}/\sqrt{p}$ at $p=10^3$, while the previous column shows conjectured values. Columns 2-4 present the effective angles $\alpha_0$, $\alpha_1$, $\alpha_2$ (those that are equal to or exceed $\pi$ are shown in gray).}
\label{tab:ck2}
\end{table}

To interpret the numerical results, we introduce the notion of ``effective angles".
%
Let us denote the angles of the triangle as $\alpha_0^{(0)}=\pi/12$, $\alpha_1^{(0)}=\pi/3$ and $\alpha_2^{(0)} = 7\pi/12$. We conjecture that the coefficient $c_0$ is given by the smallest angle $\alpha_0^{(0)}$: $c_0=\sin(\alpha_0^{(0)}/2) = \sin(\pi/24) \approx 0.1305$. The next coefficient $c_1$ turns out to be close to $c_1 \approx \sin(3\alpha_0^{(0)}/2) = \sin(\pi/8) \approx 0.3827$, i.e. one replaces the original angle $\alpha_0^{(0)}$ by $3\alpha_0^{(0)}$, by adding $2\alpha_0^{(0)}$. In other words, one now deals with effective angles $\alpha_0^{(1)}=3\pi/12$, $\alpha_1^{(1)}=\pi/3$, $\alpha_2^{(1)}=7\pi/12$, and the smallest of them determines $c_1$. Then, we apply this iterative re-adjustment of the effective angles, i.e. we identify the index $i$ of the smallest effective angle and update it by adding twice larger original angle with index $i$: 
\begin{equation}
\alpha_i^{(k+1)} = \alpha_i^{(k)} + 2\alpha_i^{(0)}.
\end{equation}
For instance, we get for $k=2$: $\alpha_0^{(2)}=5\pi/12$, $\alpha_1^{(2)}=\pi/3$, $\alpha_2^{(2)}=7\pi/12$, so that $c_2$ is determined by the smallest angle $\alpha_1^{(2)}$: $c_2\approx \sin(\alpha_1^{(2)}/2) = \sin(\pi/6) = 0.5$. However, when an effective angle exceeds $\pi$, it does not contribute to the iterative procedure anymore. The effectives angles and the conjectured values of $c_k$ are given in \Cref{tab:ck2}.
\subsection{Conjecture for polygonal domains}
The above iterative process can be generalized to arbitrary polygonal domains.
Let $a_0 =
\{\alpha_0,\ldots,\alpha_{N-1}\}$
be a sequence of all angles of a polygonal domain.
The coefficient $c_0$ is set to be $c_0 = \sin(\min\{\pi, a_0\}/2)$,
i.e., it
is determined by the smallest angle of the domain, say, $\alpha_i$.  
After
that, the $i$-th element of the sequence $a_0$ is increased by
$2\alpha_i$
to get the updated sequence $a_1$ that determines the next coefficient:
$c_1 = \sin(\min\{\pi, a_1\}/2)$. At step $k$, one sets $c_k =
\sin(\min\{\pi, a_k\}/2)$,
i.e., the smallest angle in the already constructed sequence $a_k$. If
this
angle stands on the position $j$, the $j$-th element of this sequence is
increased by $2\alpha_j$ to produce a new sequence $a_{k+1}$, and so on. Note that if there are several equal minima in the sequence of effective angles, one choose the one, for which the increment $2\alpha_i$ is the largest (and if the increments are identical, one chooses any of them).
It is clear that, after a number of steps, all angles in the sequence
will exceed $\pi$, so that all the remaining coefficients $c_k$ become equal
$1$. In the case of a regular polygon with $N$ sides and equal angles $\alpha$, the conjecture simply implies that the first $N$ coefficients are $c_k=\sin(\alpha/2)$, while the other $c_k=1$, in agreement with our earlier numerical results. 
We check the conjecture by considering a generic polygonal domain with angles $\pi/12$, $\pi/12$, $\pi/4$, $\pi/4$, $1.9064$, $2.7224$, $23\pi/12$, $23\pi/12$ (\cref{sfig:poly8}). \Cref{tab:ck4} presents the effective angles, the conjectured coefficients $c_k$ and those obtained numerically. 

Note that the smallest effective angle that determines the coefficient $c_k$, also indicates the region in which the related eigenfunctions $V_k^{(p)}$ is concentrated in the limit $p\to\infty$. In our examples, the value $p=10^3$ provided very good agreement between theoretical and numerical results but one may need to increase the value of $p$ for other domains.

\begin{table}[!htp]
\centering
\renewcommand{\arraystretch}{2}
$\begin{array}{|c|cccccccc|cc|}
\hline \boldsymbol{k}  & \boldsymbol{\alpha_0} & \bf \boldsymbol{\alpha_1} & \boldsymbol{\alpha_2} & \boldsymbol{\alpha_3} & \boldsymbol{\alpha_4} & \boldsymbol{\alpha_5} & \boldsymbol{\alpha_6} & \boldsymbol{\alpha_7} & \boldsymbol{c_k} \textbf{(conjecture)} & \boldsymbol{c_k} \textbf{(numeric)}  \\
\hline 
0 & \cfrac{\pi}{12} & \cfrac{\pi}{12} & \cfrac{\pi}{4} & \cfrac{\pi}{4} & 1.9064 & 2.7224 & \textcolor{gray}{\cfrac{23\pi}{12}} & \textcolor{gray}{\cfrac{23\pi}{12}} & \sin\left(\cfrac{1}{2}\cdot\cfrac{\pi}{12}\right)=0.1305 &  0.1306 \\
1 & \cfrac{3\pi}{12} & \cfrac{\pi}{12} & \cfrac{\pi}{4} & \cfrac{\pi}{4} & 1.9064 & 2.7224 & \textcolor{gray}{\cfrac{23\pi}{12}} & \textcolor{gray}{\cfrac{23\pi}{12}} & \sin\left(\cfrac{1}{2}\cdot\cfrac{\pi}{12}\right)=0.1305 &  0.1306 \\
2 & \cfrac{3\pi}{12} & \cfrac{3\pi}{12} & \cfrac{\pi}{4} & \cfrac{\pi}{4} & 1.9064 & 2.7224 & \textcolor{gray}{\cfrac{23\pi}{12}} & \textcolor{gray}{\cfrac{23\pi}{12}} & \sin\left(\cfrac{1}{2}\cdot\cfrac{\pi}{4}\right)=0.3827 &  0.3828 \\
3 & \cfrac{3\pi}{12} & \cfrac{3\pi}{12} & \cfrac{3\pi}{4} & \cfrac{\pi}{4} & 1.9064 & 2.7224 & \textcolor{gray}{\cfrac{23\pi}{12}} & \textcolor{gray}{\cfrac{23\pi}{12}} & \sin\left(\cfrac{1}{2}\cdot\cfrac{\pi}{4}\right)=0.3827 &  0.3828 \\
4 & \cfrac{3\pi}{12} & \cfrac{3\pi}{12} & \cfrac{3\pi}{4} & \cfrac{3\pi}{4} & 1.9064 & 2.7224 & \textcolor{gray}{\cfrac{23\pi}{12}} & \textcolor{gray}{\cfrac{23\pi}{12}} & \sin\left(\cfrac{1}{2}\cdot\cfrac{3\pi}{12}\right)=0.3827 &  0.3838 \\
5 & \cfrac{5\pi}{12} & \cfrac{3\pi}{12} & \cfrac{3\pi}{4} & \cfrac{3\pi}{4} & 1.9064 & 2.7224 & \textcolor{gray}{\cfrac{23\pi}{12}} & \textcolor{gray}{\cfrac{23\pi}{12}} & \sin\left(\cfrac{1}{2}\cdot\cfrac{3\pi}{12}\right)=0.3827 &  0.3838 \\
6 & \cfrac{5\pi}{12} & \cfrac{5\pi}{12} & \cfrac{3\pi}{4} & \cfrac{3\pi}{4} & 1.9064 & 2.7224 & \textcolor{gray}{\cfrac{23\pi}{12}} & \textcolor{gray}{\cfrac{23\pi}{12}} & \sin\left(\cfrac{1}{2}\cdot\cfrac{5\pi}{12}\right)=0.6088 &  0.6113 \\
7 & \cfrac{7\pi}{12} & \cfrac{5\pi}{12} & \cfrac{3\pi}{4} & \cfrac{3\pi}{4} & 1.9064 & 2.7224 & \textcolor{gray}{\cfrac{23\pi}{12}} & \textcolor{gray}{\cfrac{23\pi}{12}} & \sin\left(\cfrac{1}{2}\cdot\cfrac{5\pi}{12}\right)=0.6088 &  0.6113 \\
8 & \cfrac{7\pi}{12} & \cfrac{7\pi}{12} & \cfrac{3\pi}{4} & \cfrac{3\pi}{4} & 1.9064 & 2.7224 & \textcolor{gray}{\cfrac{23\pi}{12}} & \textcolor{gray}{\cfrac{23\pi}{12}} & \sin\left(\cfrac{1}{2}\cdot\cfrac{7\pi}{12}\right)=0.7934 &  0.7967 \\
9 & \cfrac{9\pi}{12} & \cfrac{7\pi}{12} & \cfrac{3\pi}{4} & \cfrac{3\pi}{4} & 1.9064 & 2.7224 & \textcolor{gray}{\cfrac{23\pi}{12}} & \textcolor{gray}{\cfrac{23\pi}{12}} & \sin\left(\cfrac{1}{2}\cdot\cfrac{7\pi}{12}\right)=0.7934 &  0.7967 \\
10 & \cfrac{9\pi}{12} & \cfrac{9\pi}{12} & \cfrac{3\pi}{4} & \cfrac{3\pi}{4} & 1.9064 & 2.7224 & \textcolor{gray}{\cfrac{23\pi}{12}} & \textcolor{gray}{\cfrac{23\pi}{12}} & \sin\left(\cfrac{1.9064}{2}\right)=0.8153 &  0.8159 \\
11 & \cfrac{9\pi}{12} & \cfrac{9\pi}{12} & \cfrac{3\pi}{4} & \cfrac{3\pi}{4} & \textcolor{gray}{5.7192} & 2.7224 & \textcolor{gray}{\cfrac{23\pi}{12}} & \textcolor{gray}{\cfrac{23\pi}{12}} & \sin\left(\cfrac{1}{2}\cdot\cfrac{3\pi}{4}\right)=0.9239 &  0.9257 \\
12 & \cfrac{9\pi}{12} & \cfrac{9\pi}{12} & \cfrac{3\pi}{4} & \textcolor{gray}{\cfrac{9\pi}{4}} & \textcolor{gray}{5.7192} & 2.7224 & \textcolor{gray}{\cfrac{23\pi}{12}} & \textcolor{gray}{\cfrac{23\pi}{12}} & \sin\left(\cfrac{1}{2}\cdot\cfrac{3\pi}{4}\right)=0.9239 &  0.9257 \\
13 & \cfrac{9\pi}{12} & \cfrac{9\pi}{12} & \textcolor{gray}{\cfrac{9\pi}{4}} & \textcolor{gray}{\cfrac{9\pi}{4}} & \textcolor{gray}{5.7192} & 2.7224 & \textcolor{gray}{\cfrac{23\pi}{12}} & \textcolor{gray}{\cfrac{23\pi}{12}} & \sin\left(\cfrac{1}{2}\cdot\cfrac{9\pi}{12}\right)=0.9239 &  0.9285 \\
14 & \cfrac{11\pi}{12} & \cfrac{9\pi}{12} & \textcolor{gray}{\cfrac{9\pi}{4}} & \textcolor{gray}{\cfrac{9\pi}{4}} & \textcolor{gray}{5.7192} & 2.7224 & \textcolor{gray}{\cfrac{23\pi}{12}} & \textcolor{gray}{\cfrac{23\pi}{12}} & \sin\left(\cfrac{1}{2}\cdot\cfrac{9\pi}{12}\right)=0.9239 &  0.9285 \\
15 & \cfrac{11\pi}{12} & \cfrac{11\pi}{12} & \textcolor{gray}{\cfrac{9\pi}{4}} & \textcolor{gray}{\cfrac{9\pi}{4}} & \textcolor{gray}{5.7192} & 2.7224 & \textcolor{gray}{\cfrac{23\pi}{12}} & \textcolor{gray}{\cfrac{23\pi}{12}} & \sin\left(\cfrac{2.7224}{2}\right)=0.9781 & 0.9789  \\
16 & \cfrac{11\pi}{12} & \cfrac{11\pi}{12} & \textcolor{gray}{\cfrac{9\pi}{4}} & \textcolor{gray}{\cfrac{9\pi}{4}} & \textcolor{gray}{5.7192} & \textcolor{gray}{8.1672} & \textcolor{gray}{\cfrac{23\pi}{12}} & \textcolor{gray}{\cfrac{23\pi}{12}} & \sin\left(\cfrac{1}{2}\cdot\cfrac{11\pi}{12}\right)=0.9914 &  0.9946 \\
17 & \cfrac{11\pi}{12} & \textcolor{gray}{\cfrac{13\pi}{12}} & \textcolor{gray}{\cfrac{9\pi}{4}} & \textcolor{gray}{\cfrac{9\pi}{4}} & \textcolor{gray}{5.7192} & \textcolor{gray}{8.1672} & \textcolor{gray}{\cfrac{23\pi}{12}} & \textcolor{gray}{\cfrac{23\pi}{12}} & \sin\left(\cfrac{1}{2}\cdot\cfrac{11\pi}{12}\right)=0.9914 &  0.9964 \\
\hline
\end{array}$
\caption{First 18 coefficients $c_k$ for the polygon shown in \cref{sfig:poly8} with angles $\pi/12$, $\pi/12$, $\pi/4$, $\pi/4$, $1.9064$, $2.7224$, $23\pi/12$, $23\pi/12$. The last column presents the values $c_k$ estimated from the ratio $\mu_k^{(p)}/\sqrt{p}$ at $p=10^3$, while the previous column shows conjectured values. Columns 2-9 present the effective angles (those that exceed $\pi$ are shown in gray).}
\label{tab:ck4}
\end{table}

\section{Coefficients of spectral expansions}\label{sec:trunc}
In this section, we investigate how the coefficients $A_k^{(p)}$ defined by \cref{eq:Ak} depend on $p$ and $k$ for various planar domains. 
As mentioned in \cref{sec:intro}, $A_k^{(0)}=\delta_{k,0}$ due to orthogonality of eigenfunctions $v_k^{(0)}$ to $v_0^{(0)}=1/\sqrt{|\pa|}$. However, for $p>0$, the eigenfunction $v_0^{(p)}$ is in general not constant so that $A_k^{(p)}$ may be non zero for $k>0$. In the case of a disk, the rotational symmetry implies that $v_k^{(p)}$ do not depend on $p$, so that $A_k^{(p)}=\delta_{k,0}$ for any $p\geq0$. As a consequence, spectral expansions involving $A_k^{(p)}$ as coefficients are reduced to a single term \cite{grebenkov2020paradigm}. One may wonder how the domain shape can alter this behavior.

One can see the coefficient $A_k^{(p)}$ as the coefficients in the expansion of a constant function $1/\sqrt{|\pa|}$ over the orthogonal basis $\{v_k^{(p)}\}$ of $L^2(\pa)$. As a consequence, one has 
\begin{equation}
\sum_{k=0}^{\infty} |A_k^{(p)}|^2=1.
\end{equation}
In particular,  $|A_k^{(p)}|^2$ can be interpreted as relative weights of different eigenfunctions $v_k^{(p)}$ in some spectral expansions. The convergence of this series implies that $|A_k^{(p)}|$ decreases with $k$, i.e.
\begin{equation}\label{eq:asymp-Ak}
\lim\limits_{k\to\infty} A_k^{(p)} = 0.
\end{equation}
We inspect the relative contributions of $A_k^{(p)}$ for various domains.
\subsection{Ellipses and rectangles}
\begin{figure}[!ht]
\centering
\subfloat[
\label{fig:Ak_ellip2}]{%
\includegraphics[width=0.45\linewidth]{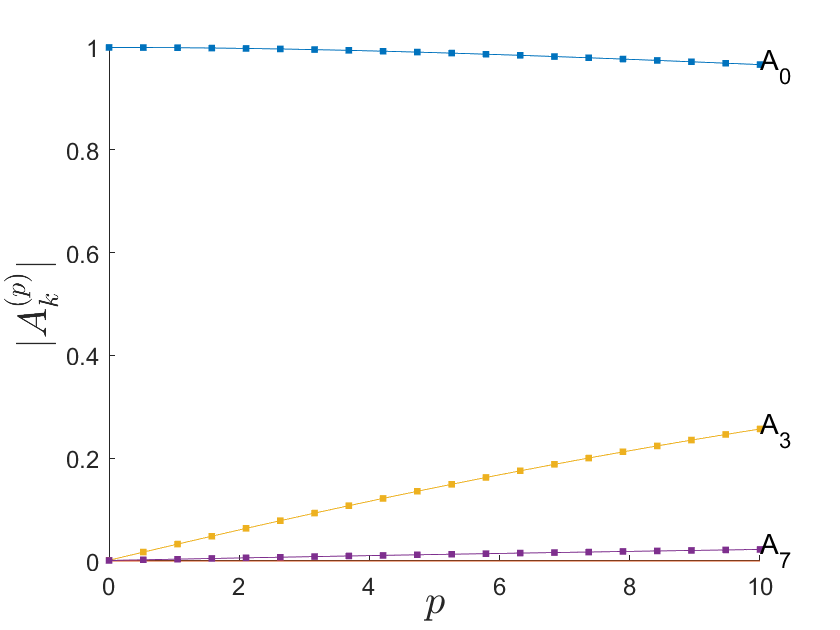}}
\quad
\subfloat[
\label{fig:Ak_ellip10}]{\includegraphics[width=0.45\linewidth]{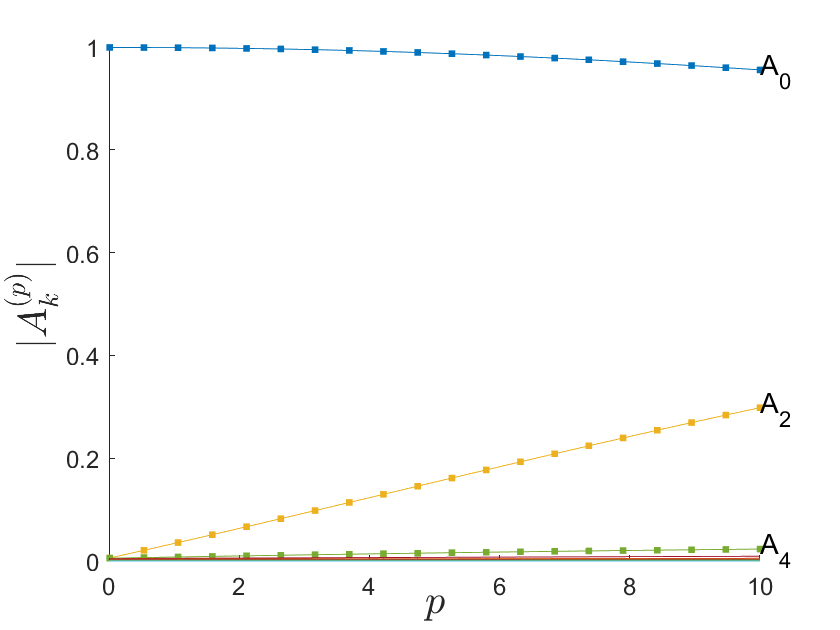}}
\caption{First 20 coefficients $|A_k^{(p)}|$ for (a) an ellipse of semiaxes 1/2 and 1 (aspect ratio 2), (b) an ellipse of semiaxes 1/10 and 1 (aspect ratio 10).}
\label{fig:Ak_ellip}
\end{figure}

In order to break the rotational invariance of a disk, we start again with ellipses of aspect ratio $2$ and $10$.
\Cref{fig:Ak_ellip} shows the first 20 coefficients $|A_k^{(p)}|$ for both ellipses. We observe that $A_k^{(p)} \approx 0$ a wide range of $k$ and $p$, except for some coefficients, namely, $A_0^{(p)}$, $A_3^{(p)}$, $A_7^{(p)}$ and $A_{11}^{(p)}$ for the ellipse of aspect ratio $2$; and $A_0^{(p)}$, $A_2^{(p)}$ and $A_4^{(p)}$ for the ellipse of aspect ratio $10$.
Despite the broken rotational invariance of $\Omega$, these results suggest that most of the coefficients $A_k^{(p)}$ vanish. This behavior is a consequence of the domain symmetry. Indeed, as an ellipse has two axes of symmetry (vertical and horizontal), the functions $V_k^{(p)}$ must be symmetric, i.e. to verify $V_k^{(p)}(-x,y) = \pm V_k^{(p)}(x,y)$ and $V_k^{(p)}(x,-y) = \pm V_k^{(p)}(x,y)$. In other words,
any $V_k^{(p)}$ satisfies these two relations with one choice of signs from 4
possible combinations. Consequently, the integration of $V_k^{(p)}$ over $\Omega$ yields $A_k^{(p)} = 0$, in 3 out of 4 cases. This is confirmed on \cref{fig:Ak_ellip2} showing the non-zero coefficients $A_0^{(p)}$, $A_3^{(p)}$, $A_7^{(p)}$, i.e. 1 out of 4. In the context of diffusion-controlled reactions, this implies that the related spectral expansions contain only a small number of contributing terms.
%

\begin{figure}[!ht]
\centering
\subfloat[\label{fig:Ak_carre}]{\includegraphics[width=0.45\linewidth]{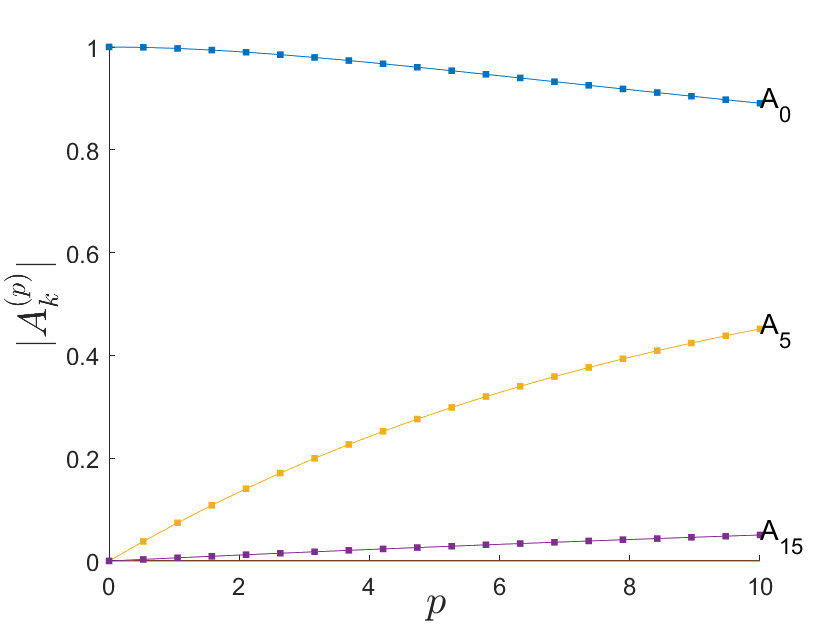}}
\quad
\subfloat[
\label{fig:Ak_rect2}]{\includegraphics[width=0.45\linewidth]{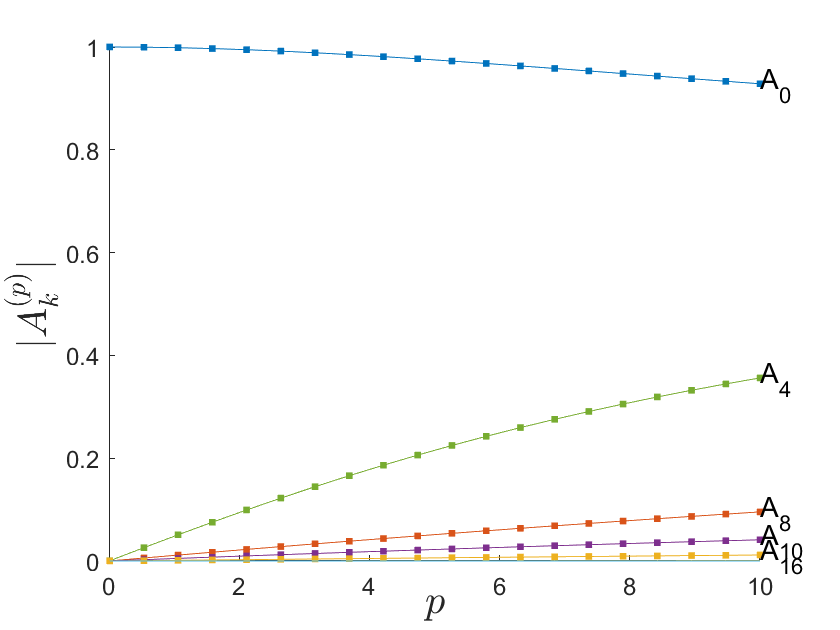}}
\caption{First 20 coefficients $|A_k^{(p)}|$ for (a) a square of sidelength 2, (b) rectangle with sides 2 and 1.}
\label{fig:Ak_rect}
\end{figure}

After revealing the role of anisotropy of the domain, we inspect to role of the smoothness of the boundary by considering rectangles (\cref{fig:Ak_rect}).
%
Despite the presence of corners, we still observe that the $A_k^{(p)}$ vanish for most $k$ and a wide range of $p$, except for some coefficients, e.g. $A_0^{(p)}$, $A_5^{(p)}$, and $A_{15}^{(p)}$ for the square, and $A_0^{(p)}$, $A_4^{(p)}$, $A_8^{(p)}$, $A_{10}^{(p)}$ and $A_{16}^{(p)}$ for the considered rectangle. This behavior can be rationalized again by the mirror symmetries of these domains.

\subsection{Prefractals}

\begin{figure}[!ht]
\centering
\subfloat[Generation 0]{\includegraphics[width=0.45\linewidth]{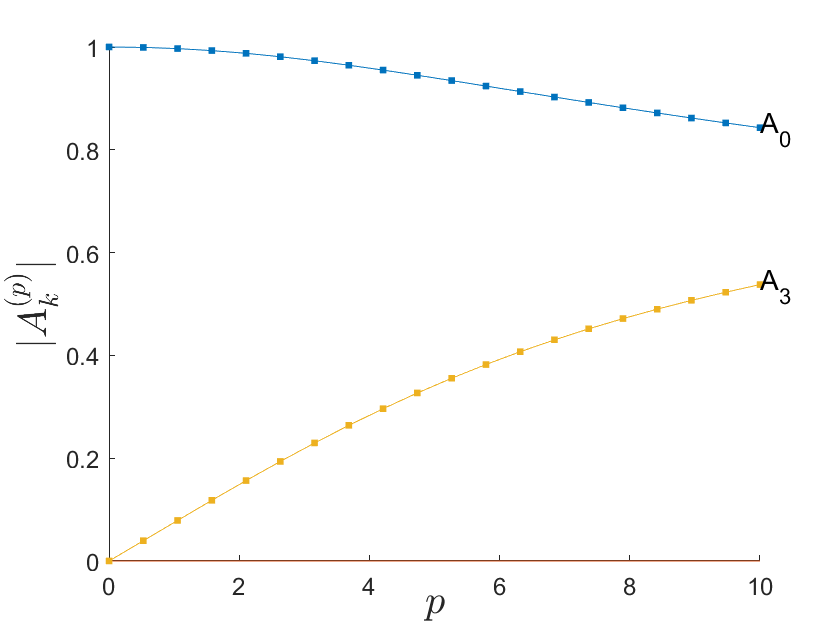}}
\hfil
\subfloat[Generation 1]{\includegraphics[width=0.45\linewidth]{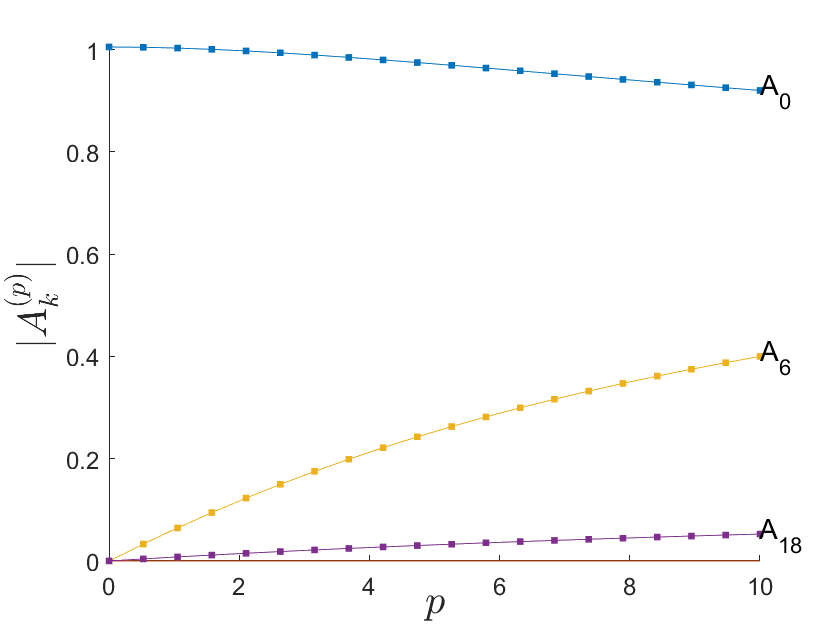}}
\quad
\subfloat[Generation 2]{\includegraphics[width=0.45\linewidth]{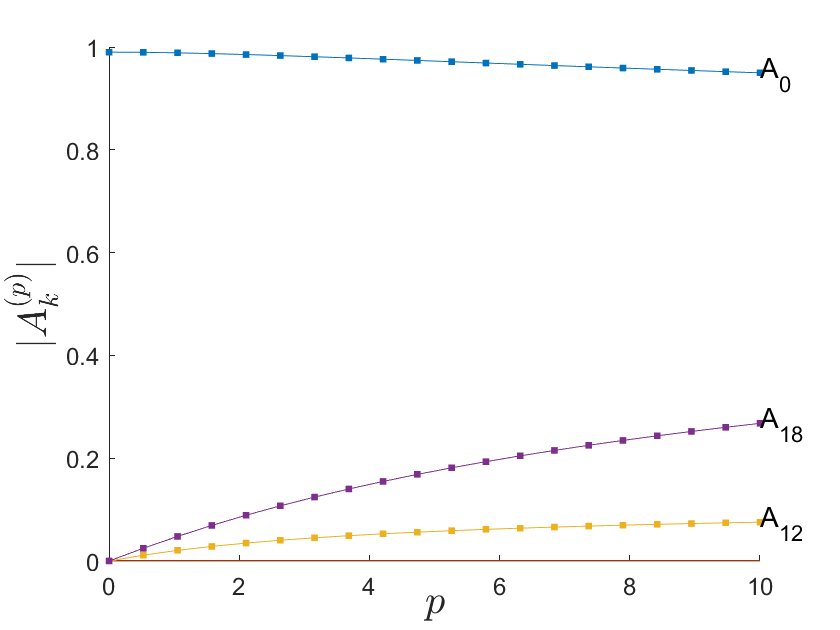}}
\subfloat[Generic triangle \label{fig:Ak_tri-pi12}]{\includegraphics[width=0.45\linewidth]{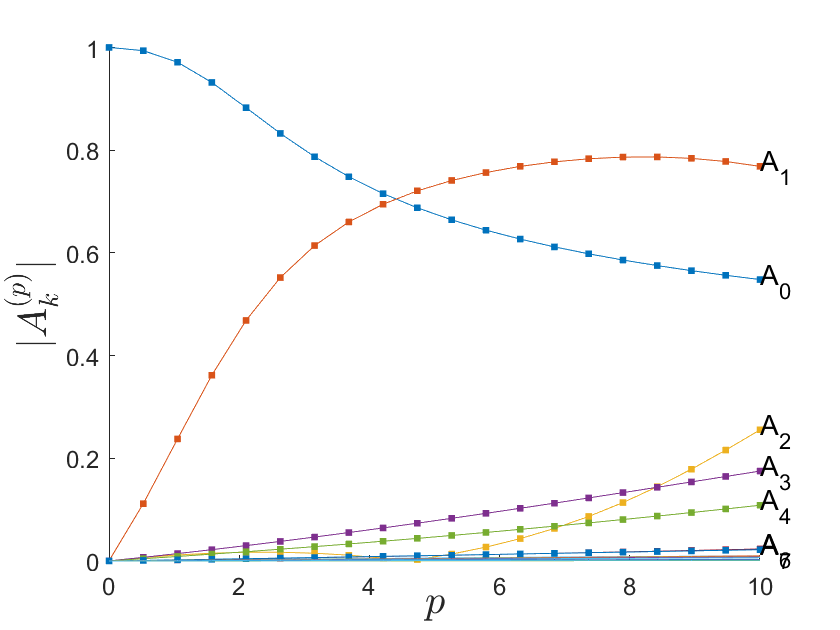}}
\caption{First 20 coefficients $|A_k^{(p)}|$ for (a), (b), (c) the first three generations of the Koch snowflake, and for (d) a triangle with sidelength 2 and angles $\pi/3$ and $\pi/12$.}
\label{fig:Ak_fractal}
\end{figure}

\Cref{fig:Ak_fractal} presents the first 20 coefficients $|A_k^{(p)}|$ for the first three generations of the Koch snowflake.
Despite the increased roughness of the prefactal boundary, we still obtain $A_k^{(p)} = 0$ for a wide range of $k$ and $p$, except for $A_0^{(p)}$ and $A_3^{(p)}$ for the equilateral triangle (generation 0); $A_0^{(p)}$, $A_6^{(p)}$ and $A_{18}^{(p)}$ for the generation 1; $A_0^{(p)}$, $A_{12}^{(p)}$ and $A_{18}^{(p)}$ for the generation 2.
As previously, cancellation of most coefficients $A_k^{(p)}$ is a consequence of symmetries, which are, however, more sophisticated than in the previous cases of ellipses and rectangles.

\subsection{Generic triangles}
To highlight the role of symmetries in previous examples, we consider a generic triangle with one sidelength 2 and two angles $\pi/3$ and $\pi/12$ shown on \cref{sfig:tri}. \Cref{fig:Ak_tri-pi12} presents the first 20 coefficients $|A_k^{(p)}|$ but their behavior is very different from the previous ones. Indeed, many coefficients $A_k^{(p)}$ are no longer close to 0. Curiously, the principal eigenfunction $V_0^{(p)}$, which provided the dominant contribution in all previous examples, gives away its leading role to the next eigenfunction $V_1^{(p)}$ for $p\geq4$.

We conclude that cancellation of many $A_k^{(p)}$ for $p>0$, which was observed for symmetric domains such as ellipses, rectangles, or even prefractal shapes, is not generic. For an arbitrary domain (e.g. a generic triangle) many terms can contribute to spectral expansions involving $A_k^{(p)}$, even thought their contributions are necessarily reduced as $k\to\infty$. As a consequence, the effect of domain geometry onto various characteristics of diffusion-controlled reactions can be much more sophisticated and versatile than one might expect from earlier theoretical studies focused on highly symmetric domains. Further numerical investigations of this effect and more rigorous characterization of the decay of $|A_k^{(p)}|$ with $k$ present an important perspective for future research.

\section{Asymptotic behavior of Steklov eigenfunctions}\label{sec:asymp2}

In this section, we analyze the decay of the Steklov eigenfunctions $V_k^{(p)}$ away from the boundary \cite{hislop2001spectral,polterovich2019nodal,galkowski2019pointwise,daude2021exponential,helffer2022semi}. In particular, Polterovich, Sher and Toth \cite{polterovich2019nodal} proved that for any bounded domain $\Omega \subseteq \mathbb{R}^2$ with a real-analytic boundary $\pa$, there exist positive constants $\eta>0$ and $B>0$ depending only on the geometry of $\Omega$, such that all eigenfunctions $V_k^{(0)}$ satisfy 

\begin{equation}\label{eq:polt}
\sqrt{|\pa|}~|V_k^{(0)}(\x)| \leq B \exp(-\eta~\mu_k^{(0)} |\x-\pa|),
\end{equation}
where $|\x-\partial \Omega|$ is the Euclidean distance between $\x$ and the boundary $\pa$ and we included the prefactor $\sqrt{|\pa|}$ to ensure that the left-hand side is dimensionless. In other words, each Steklov eigenfunction $V_k^{(0)}$ has an upper bound, which decays exponentially fast away from the boundary, and its decay rate is proportional to the corresponding eigenvalue $\mu_k^{(0)}$. As $\mu_k^{(0)}$ grows with $k$ up to infinity, the exponential decay becomes faster and more and more restrictive. In contrast, the upper bound \cref{eq:polt} is in general not much informative for any finite $k$: as $V_k^{(0)}$ is an analytic function on a bounded domain, one can always choose a large enough $B$ or small enough $\eta$ to fulfill the inequality \cref{eq:polt}. In other words, without restricting estimates on $\eta$ and $B$, the upper bound (\ref{eq:polt}) does not tell much on the behavior of a given Steklov eigenfunction. 
The upper bound (\ref{eq:polt}) was further generalized to higher dimension and $p\neq0$ by Helffer and Kachmar \cite{helffer2022semi}. They proved that for any bounded domain $\Omega \subset \mathbb{R}^d$ ($d\geq2$), with a real-analytic boundary $\pa$, any eigenfunction $V_k^{(p)}$ with sufficiently high eigenvalue $\mu_k^{(p)}$ (i.e. sufficiently high $k$)  has an upper bound that decays exponentially fast away from the boundary. In fact, for any $p>-\lambda_0^D$ (where $\lambda_0^D$ is the smallest eigenvalue of the Dirichlet Laplacian in $\Omega$) there exist constants $B>0$, $\eta>0$, $\epsilon>0$ and $k_0$ such that
\begin{equation}\label{eq:helffer}
\forall k>k_0,~\forall \x \in \Omega, \quad  |\sqrt{\pa}||V_k^{(p)}(\x)| \leq B \left(\mu_k^{(p)}\right)^{\frac{d}{2}-\frac{1}{4}} \exp \left(-\eta~ \mu_k^{(p)} \min \{\epsilon,|\x-\partial \Omega|\}\right),
\end{equation}
where we included again the prefactor $\sqrt{|\pa|}$.
Helffer and Kachmar also questioned whether the assumption of real-analytic boundary could be relaxed.

Being inspired by these fascinating results, we aim to check numerically to which extend the upper bounds \cref{eq:polt,eq:helffer} determine the decay of any given Steklov eigenfunction $V_k^{(p)}$ (with finite $k$). One may wonder whether the exponential decay holds in the whole domain, and whether the eigenvalue $\mu_k^{(p)}$ is the decay rate, i.e. whether $\eta$ is close to 1. To gain intuitive insights onto this behavior, one can first look at the explicit \cref{eq:eigelem} for the Steklov eigenfunctions $V_k^{(p)}$ for the disk of radius $R$. In the limit $p\to0$ one has for any $k=1,2,3,\ldots$
\begin{align}\label{eq:majVk}
\nonumber |V_{2k}^{(0)}(\x)| = \frac{1}{\sqrt{\pi R}}\left(\frac{|\x|}{R}\right)^{k}&=\frac{1}{\sqrt{\pi R}} \exp\left[k \ln\left(\frac{R-|\x-\pa|}{R}\right)\right]\\
&\leq \frac{1}{\sqrt{\pi R}} \exp\left(-\frac{k}{R}|\x-\pa|\right) = \frac{\exp\left(-\mu_{2k}^{(0)}|\x-\pa|\right)}{\sqrt{\pi R}},
\end{align}
where we used that $\mu_{2k}^{(0)}=k/R$; and the same applies for $V_{2k+1}^{(0)}(\x)$.
One sees how $\mu_k^{(0)}$ controls the exponential decay of $V_k^{(0)}$ away from the boundary, with $\eta=1$ and $B=1/\sqrt{\pi R}$. In the opposite limit $p\to\infty$, the asymptotic behavior of the modified Bessel functions $I_k(z)$ implies for $r\sqrt{p}\gg1$:
\begin{equation}
\begin{aligned}
|V_k^{(p)}(\x)| \lesssim \frac{1}{\sqrt{\pi R}} \frac{\sqrt{R}}{\sqrt{r}} \exp\left(-\sqrt{p} (R-r)\right) &\simeq \frac{1}{\sqrt{\pi R}} \frac{\sqrt{R}}{\sqrt{r}} \exp\left(-\mu_k^{(p)} |\x-\pa|\right) \\
&\leq \frac{p^{1/4}}{\sqrt{\pi}}\exp\left(-\mu_k^{(p)} |\x-\pa|\right),
\end{aligned}
\end{equation}
where we used the asymptotic relation (\ref{eq:asymp-mu_circle}). 
In turn, if $r$ is so small that $r\sqrt{p} \ll 1$ (but still $R\sqrt{p}\gg 1$) one uses 
$I_k(r\sqrt{p}) \simeq \cfrac{(r\sqrt{p})^{k}}{k!}$ to get
\begin{equation}
|V_{2k}^{(p)}(\x)| \ll \frac{p^{1/4}}{k!}\exp\left(-\sqrt{p}R\right) \lesssim \frac{p^{1/4}}{k!} \exp\left(-\mu_{2k}^{(p)} |\x-\pa|\right).
\end{equation}
In both limits, we observe an exponential decay of $V_k^{(p)}$ away from the boundary, which is controlled by $\mu_k^{(p)}$, and $\eta$ is close to $1$ (even though the above upper bounds suggest $\eta=1$, our derivation involved some asymptotic relations, which may require having $\eta$ slightly smaller than 1). In Appendix \ref{secap:loc-rect}, we analyze the explicit representation of the Steklov eigenfunctions for a rectangle and come to the same conclusion for $p=0$. In sharp contrast to the disk, the boundary of the rectangle is not real-analytic due to the corners. Nevertheless, the upper bound \cref{eq:polt} holds and the eigenfunctions exhibit exponential decay controlled by $\mu_k^{(0)}$, with $\eta$ close to 1.

Let us give another argument in favor of the particular choice $\eta=1$. Let $\x_0$ be a boundary point, and $\x_{\delta} = \x_0 - \delta n_{\x_0} \in \Omega$ be a bulk point at small distance $\delta$ from $\pa$, where $n_{\x_0}$ is the normal unit vector to $\pa$ at $\x_0$ oriented outward the domain. If we assume that $V_k^{(p)}(\x_0) >0$, then the Steklov boundary condition can be written as 
\begin{equation}\label{eq:lnSteklov}
\partial_n \ln\left(V_k^{(p)}(\x_0)\right) = \mu_k^{(p)} \quad (\x_0 \in \pa)
\end{equation}
(if $V_k^{(p)}(\x_0)<0$, one can replace $V_k^{(p)}(\x_0)$ by $-V_k^{(p)}(\x_0)$ that yields \\ $\partial_n \left.\ln\left(|V_k^{(p)}(\x_0)|\right)\right|_{\pa} = \mu_k^{(p)}$, with no changes in the argument; in turn, we ignore here the specific points at which $V_k^{(p)}(\x_0)=0$). Since $V_k^{(p)}$ is analytic in $\Omega$ and $V_k^{(p)}(\x_0) >0$, it is also positive in a small vicinity of  the boundary, so that 
\begin{equation}
\ln\left(V_k^{(p)}(\x_{\delta})\right) = \ln\left(V_k^{(p)}(\x_0)\right) - \mu_k^{(p)}\delta + O(\delta^2),
\end{equation}
to be consistent with \cref{eq:lnSteklov}. As a consequence, the Steklov eigenfunction $V_k^{(p)}$ exhibits an exponential decay near the boundary
\begin{equation}\label{eq:expdecay}
V_k^{(p)}(\x_{\delta}) \approx V_k^{(p)}(\x_0) \exp(-\mu_k^{(p)} \delta),
\end{equation}
where $\delta = |\x_{\delta}-\pa|$ is the distance to the boundary. One sees that the value $\eta=1$ naturally comes from the Steklov condition. The fundamental question is whether the exponential behavior (\ref{eq:expdecay}) holds \textit{approximately} far from the boundary, when the above argument is not applicable.

\begin{figure}[!ht]
\centering
\subfloat[]{\includegraphics[width=0.33\linewidth]{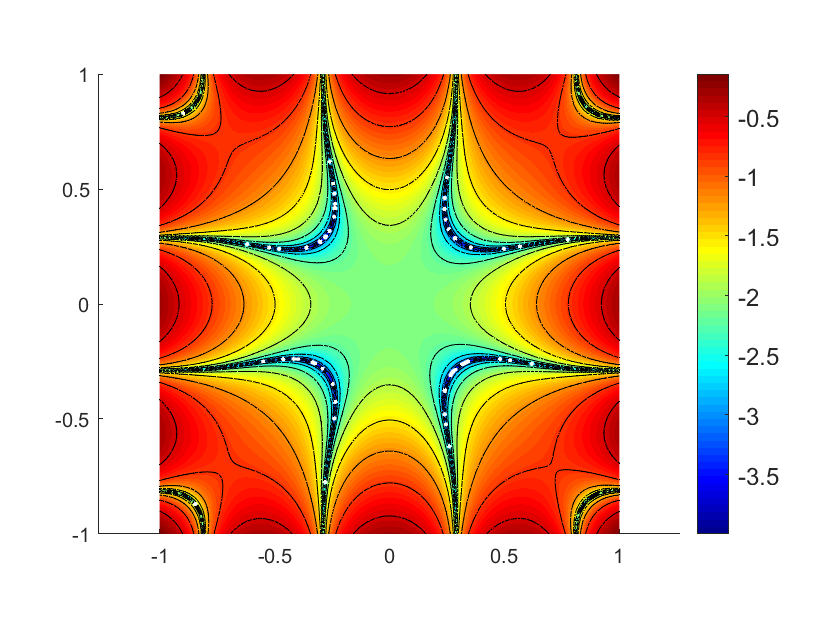}}
\subfloat[]{\includegraphics[width=0.33\linewidth]{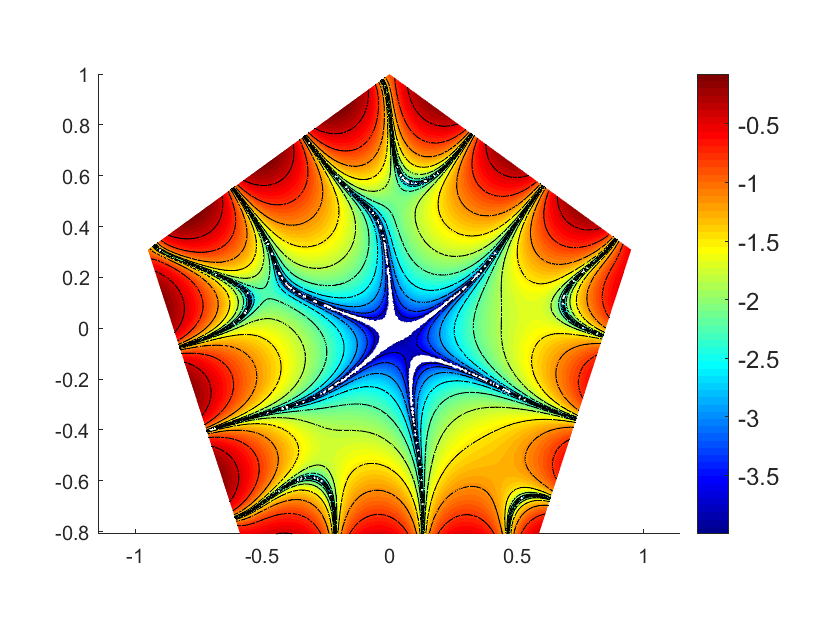}}
\subfloat[]{\includegraphics[width=0.33\linewidth]{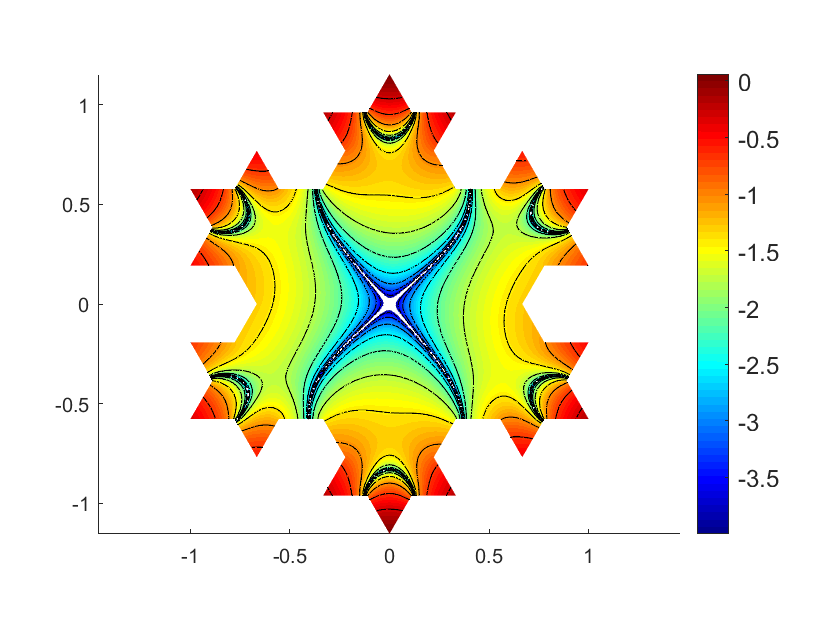}}
\hfil
\subfloat[]{\includegraphics[width=0.33\linewidth]{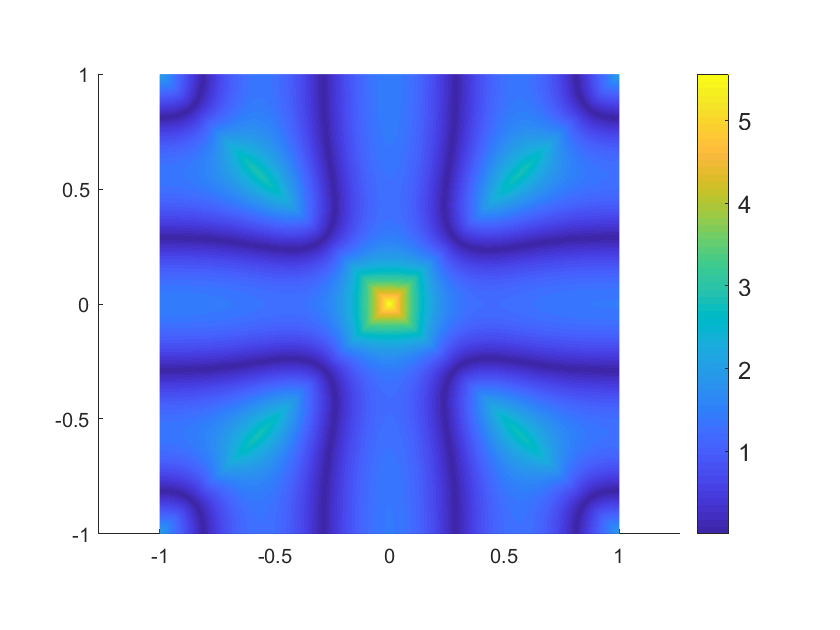}}
\subfloat[]{\includegraphics[width=0.33\linewidth]{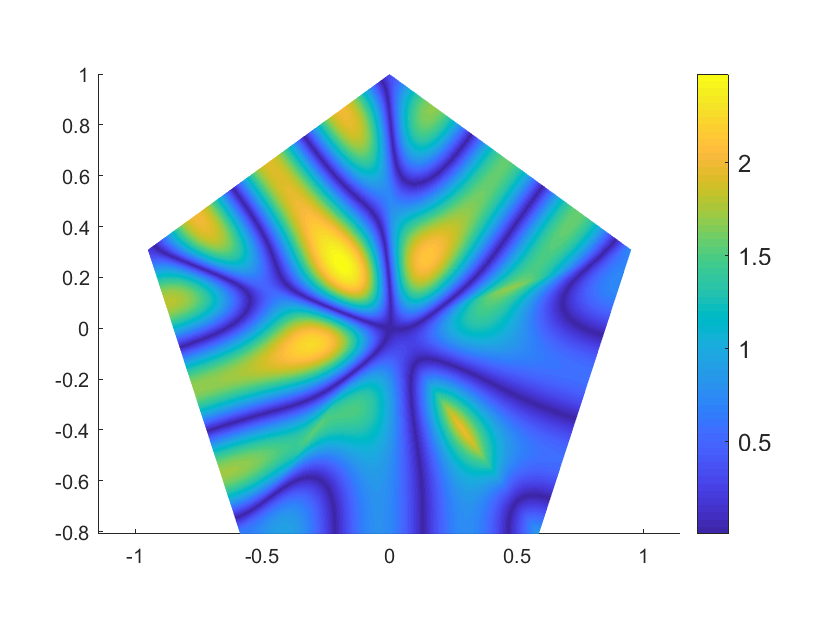}}
\subfloat[]{\includegraphics[width=0.33\linewidth]{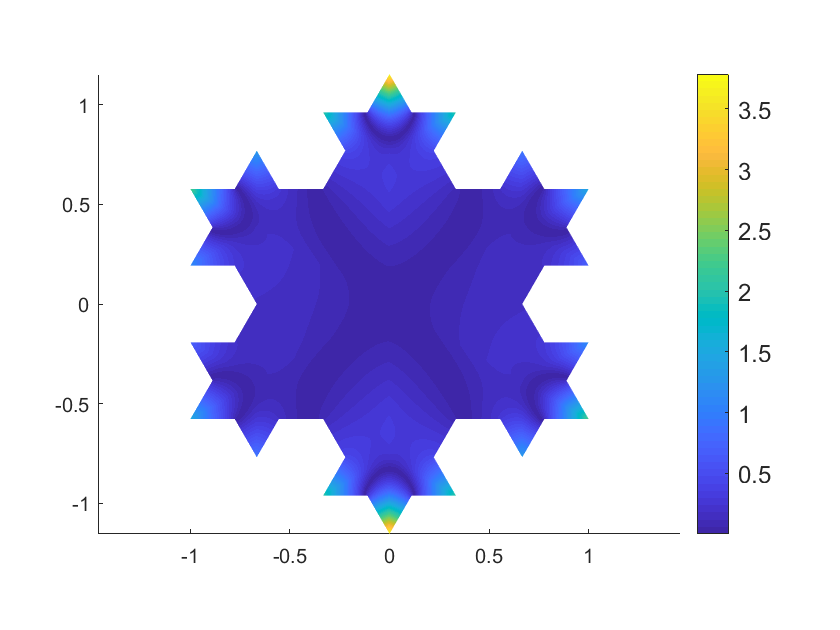}}
\caption{(Top row) Localization and decay of $V_{15}^{(0)}$ away from the boundary for (a) a square, (b) a pentagon, (c) generation 2 of the Koch snowflake. Colormaps illustrate the behavior of $\lg(|V_{15}^{(0)}|)$, where $\lg(z) = \ln(z)/\ln(10)$ is the decimal logarithm; white regions correspond to the values $|V_{15}^{(0)}|<10^{-4}$. The related eigenvalues are (a) $\mu_{15}^{(0)}\approx 5.50$, (b) $\mu_{15}^{(0)}\approx 9.00$ and (c) $\mu_{15}^{(0)}\approx 2.75$. (Bottom row) The function $B_{15}(\x)$ from \cref{eq:Bk}. \\
}
\label{fig:exploc_poly}
\end{figure}

We address this question numerically for $p=0$ and check the exponential decay of $V_k^{(0)}$ away from the boundary for various planar domains, even when $\pa$ is not real-analytic. Most importantly, we discuss whether $\eta$ is close to 1 in general. The top row of \cref{fig:exploc_poly} presents the log-scaled eigenfunction $V_{15}^{(0)}$ for three polygonal domains. It confirms the expected exponential decay away from the boundary for the square (see Appendix \ref{secap:loc-rect}), but also shows the exponential decay of $V_{15}^{(0)}$ for a pentagon and the second generation of the Koch snowflake. In these cases, the presence of corners does not seem to affect the localization near the smooth parts of the polygonal boundary.  For a more systematic insight, we define 
\begin{equation}\label{eq:Bk}
B_k(\x) =  \left| \sqrt{|\pa|}~V_k^{(0)}(\x) \exp\left(\mu_k^{(0)} |\x-\pa|\right)\right|,
\end{equation}
as if we explicitly set $\eta=1$ in \cref{eq:polt}. If $B_k(\x)$ was constant, the eigenfunction $V_k^{(0)}(\x)$ would exhibit the exponential decay away from the boundary with the rate $\mu_k^{(0)}$. As a consequence, variations of $B_k(\x)$ and, in particular, its high values  can indicate regions where the exponential decay $\exp\left(-\mu_k^{(0)} |\x - \partial\Omega|\right)$ does not hold. We note that, even for a disk, oscillations of $V_k^{(0)}(\x)$ on the boundary and inside the domain result in variations of $B_k(\x)$. 
The bottom row of \cref{fig:exploc_poly} presents the function $B_{15}(\x)$ for the three polygonal domains. For instance, we get  $\max\limits_{\x\in\Omega}B_{15}(\x) \approx 5.6$, $\max\limits_{\x\in\Omega}B_{15}(\x) \approx 2.6$ and $\max\limits_{\x\in\Omega}B_{15}(\x) \approx 3.8$ for the square, the pentagon and the Koch snowflake respectively. Interestingly, the maximal deviation from the upper bound can be found in the center (for a square), in the central part (for a pentagon), or near the boundary (for a Koch snowflake). 
However, for the pentagon and the Koch snowflake, we also observe that there exist indices $k$ for which the maximum value of $B_k(\x)$ becomes much larger, suggesting deviations of the related eigenfunctions from the exponential decay with $\eta$ close to 1.

\begin{figure}[!ht]
\centering
\subfloat[]{\includegraphics[width=0.33\linewidth]{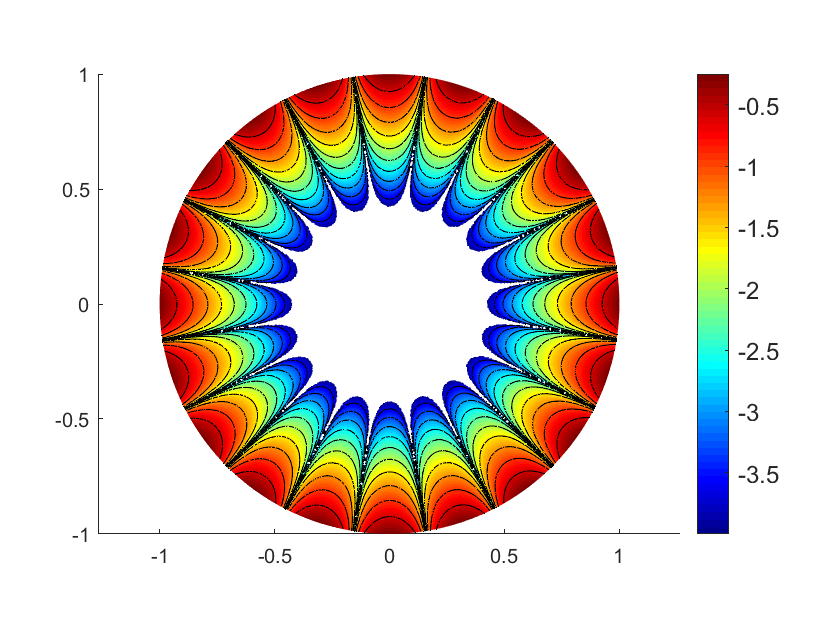}}
\subfloat[]{\includegraphics[width=0.33\linewidth]{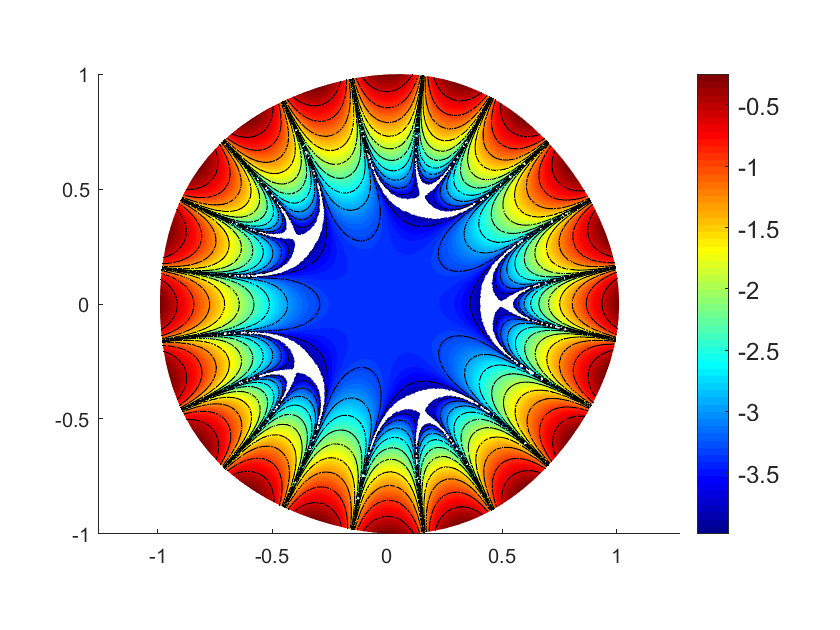}}
\subfloat[]{\includegraphics[width=0.33\linewidth]{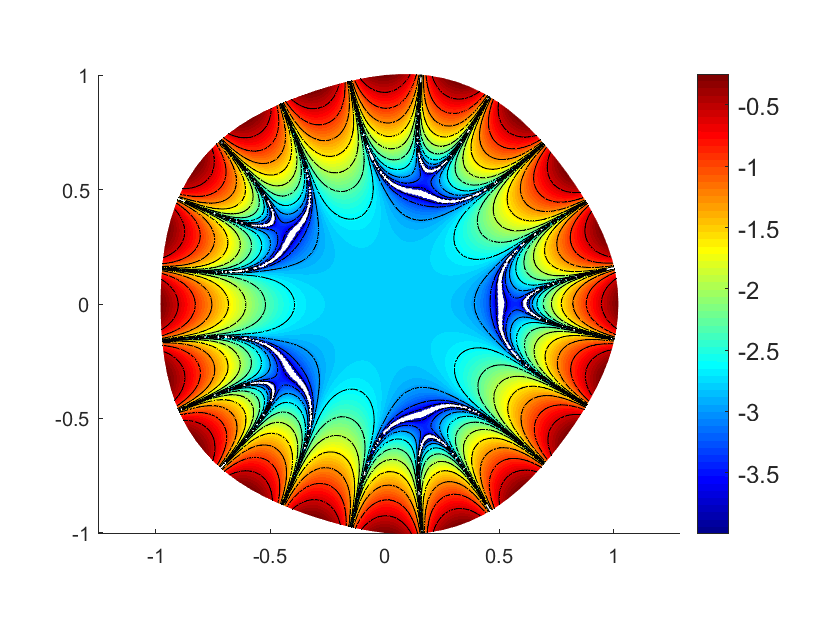}}
\hfil
\subfloat[]{\includegraphics[width=0.33\linewidth]{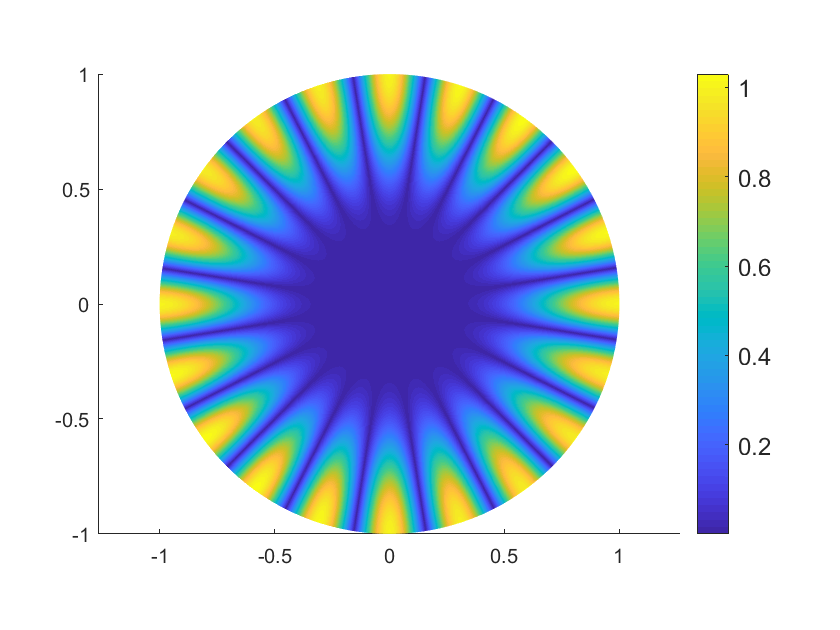}}
\subfloat[]{\includegraphics[width=0.33\linewidth]{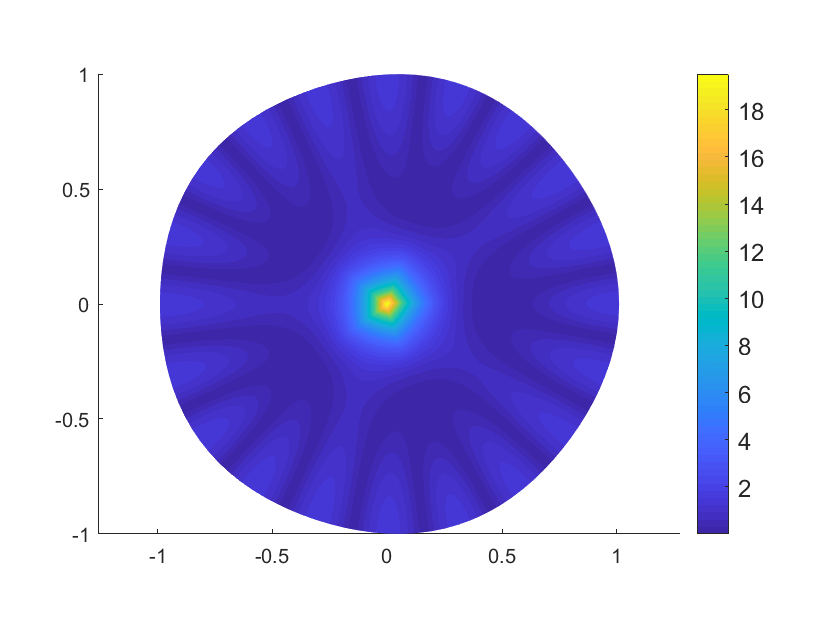}}
\subfloat[]{\includegraphics[width=0.33\linewidth]{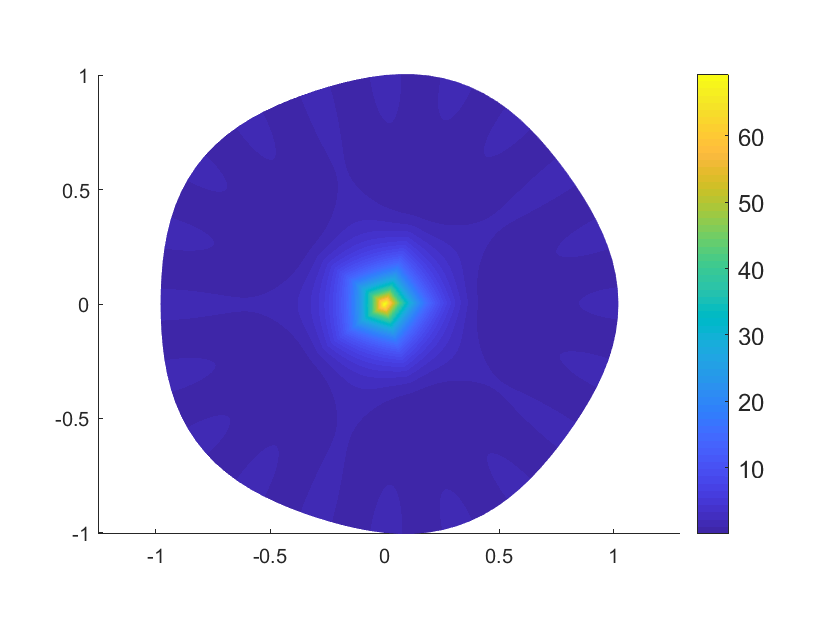}}
\caption{(Top row) Localization and decay of $V_{20}^{(0)}$ away from the boundary for (a) a disk, and two slightly deformed disks defined by \cref{eq:radius} with (b) $\gamma=0.01$, and (c) $\gamma=0.02$. Colormaps illustrate the behavior of $\lg(|V_{20}^{(0)}|)$, where $\lg(z) = \ln(z)/\ln(10)$ is the decimal logarithm; white regions correspond to the values $|V_{20}^{(0)}|<10^{-4}$. The related eigenvalues are (a) $\mu_{20}^{(0)}= 10.00$, (b) $\mu_{20}^{(0)}\approx 9.99$ and (c) $\mu_{20}^{(0)}\approx 9.98$.
(Bottom row) The function $B_{20}(\x)$ from \cref{eq:Bk}.  \\
}
\label{fig:exploc_smooth}
\end{figure}

To better understand the origin of this behavior, we go back to real-analytic boundaries and  focus on the disk and its minor deformations. Using polar coordinates $(r,\theta)$, we define the boundary of the deformed disks  $\Omega_{\gamma} = \{ (r,\theta)|~ r<\rho(\theta)\}$ by setting
\begin{equation}\label{eq:radius}
\rho(\theta) = 1 + \gamma \cos(5 \theta), \quad 0\leq \theta \leq 2\pi.
\end{equation}
\Cref{fig:exploc_smooth} presents the localization and decay of $V_{20}^{(0)}$ away from the boundary, as well as the function  $B_{20}(\x)$ for a disk, and the two slightly deformed disks with $\gamma=0.01$ and $\gamma=0.02$. As already observed in \cite{bruno2020domains}, a very small perturbation of the disk may result in significant changes of the structure of nodal lines and thus in the interior behavior of Steklov eigenfunctions. Indeed, a very small perturbation results in high values of $B_k(\x)$ at the center of the domain for some specific indices $k$, e.g. for $k=20$. 
These high values of $B_k(\x)$ suggest that the related eigenfunction decreases slower inside the domain than near the boundary; in other words, the upper bound (\ref{eq:polt}) with $\eta\approx1$ may not be valid for the whole domain. Curiously, we do not retrieve such behavior for ellipses (see Fig. 10), which could also be seen as deformations of a disk.  This  can be a consequence of symmetries.

\begin{figure}[!ht]
\centering
\subfloat[]{\includegraphics[width=0.33\linewidth]{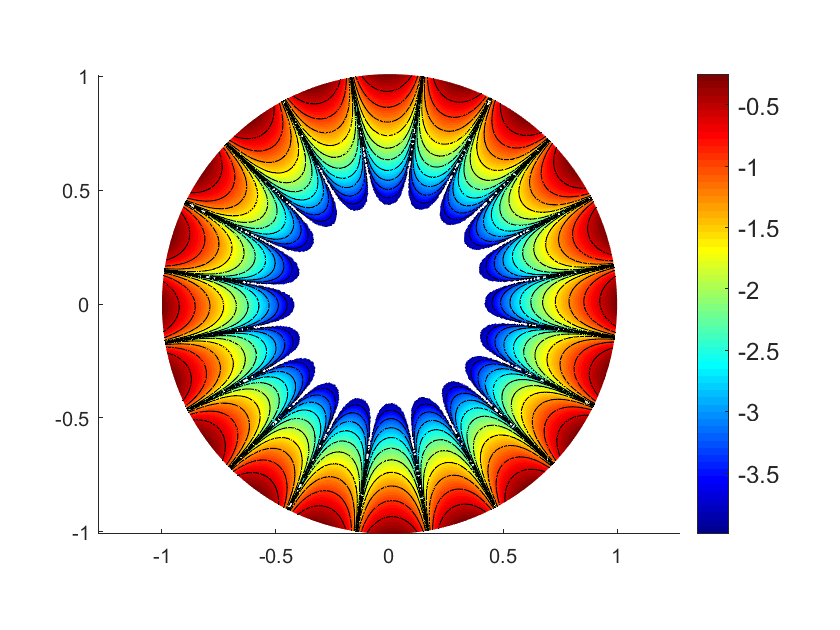}}
\subfloat[]{\includegraphics[width=0.33\linewidth]{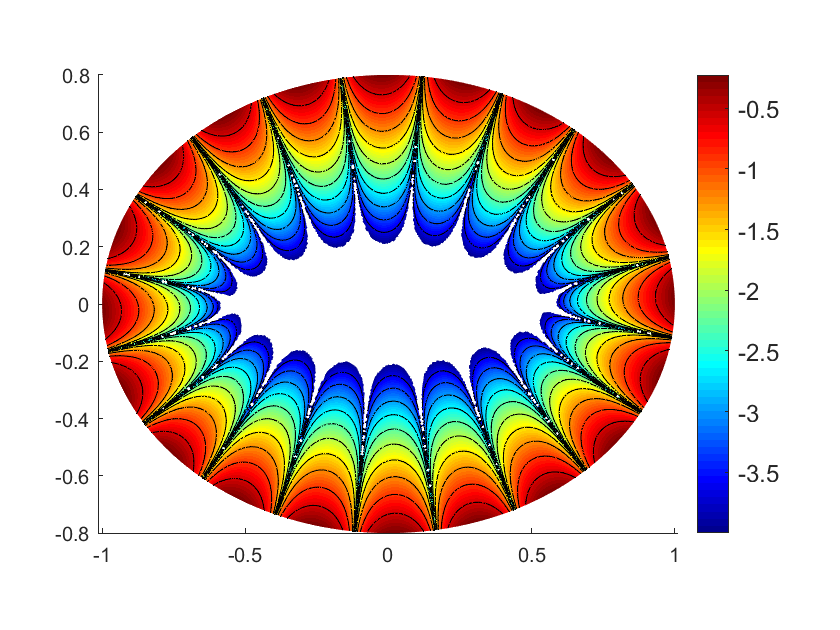}}
\subfloat[]{\includegraphics[width=0.33\linewidth]{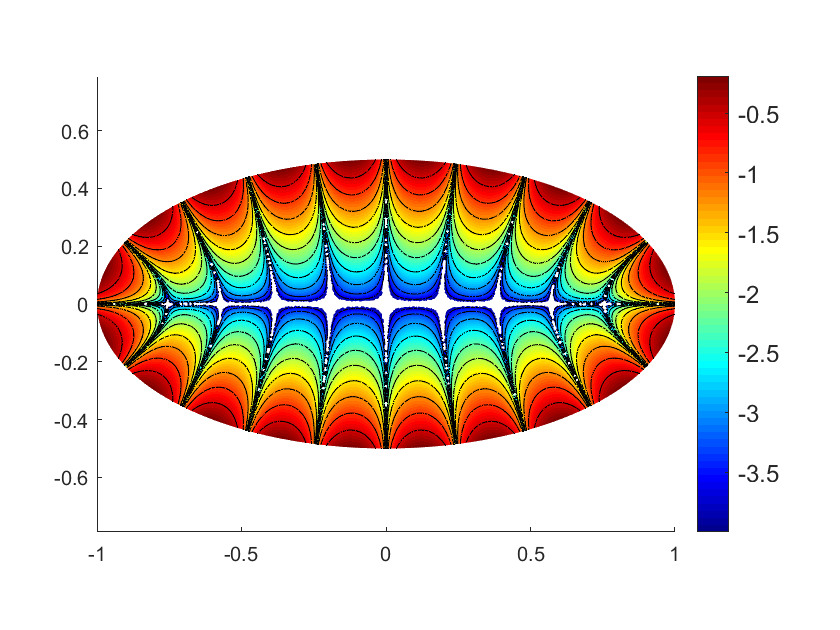}}
\hfill
\subfloat[]{\includegraphics[width=0.33\linewidth]{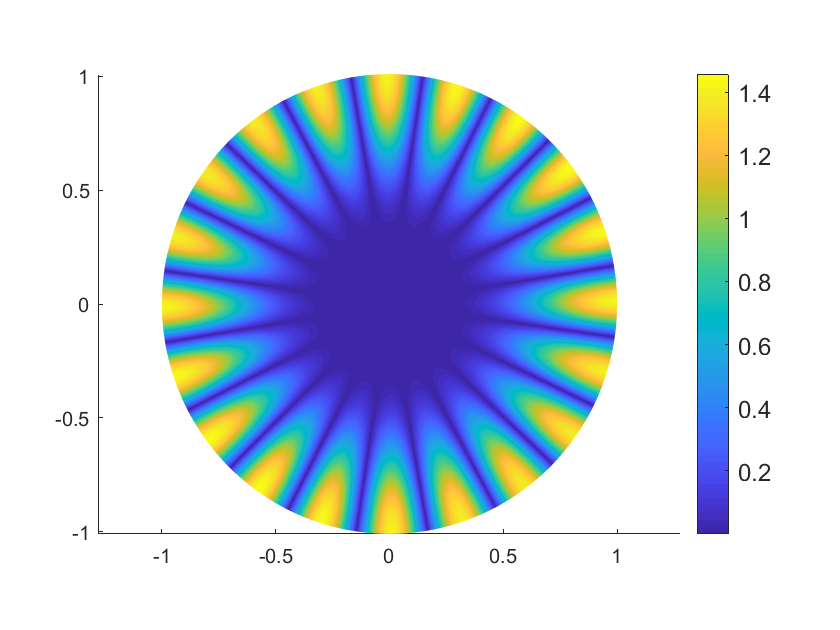}}
\subfloat[]{\includegraphics[width=0.33\linewidth]{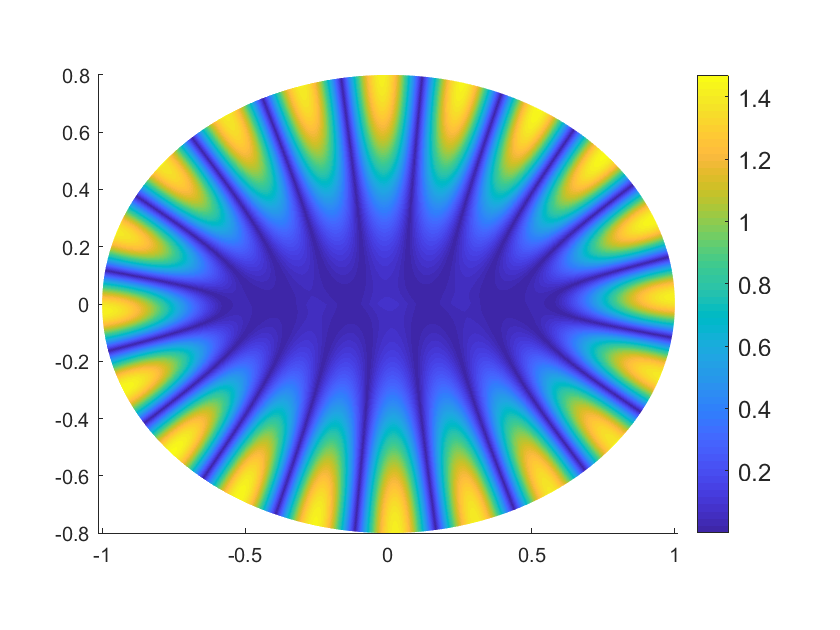}}
\subfloat[]{\includegraphics[width=0.33\linewidth]{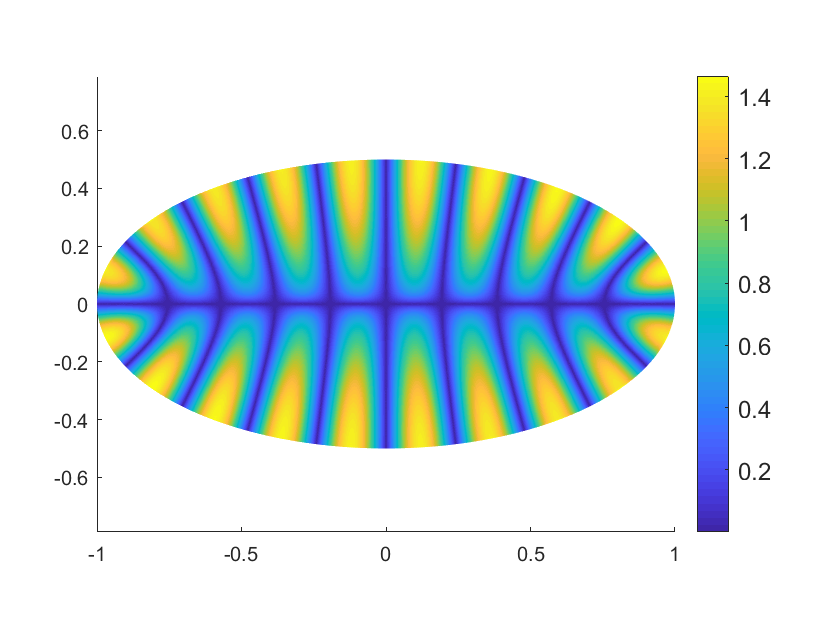}}
\caption{(Top row) Localization and decay of $V_{20}^{(0)}$ away from the boundary for (a) an ellipse with semiaxes $a=1$, $b=1.01$ (b) an ellipse with semiaxes $a=1$, $b=0.8$ , and (c) an ellipse with semiaxes $a=1$, $b=0.5$. Colormaps illustrate the behavior of $\lg(|V_{20}^{(0)}|)$. (Bottom row) The function $B_{20}(\x)$ from \cref{eq:Bk}. The related eigenvalues are (a) $\mu_{20}^{(0)}\approx 9.95$, (b) $\mu_{20}^{(0)}\approx 11.08$ and (c) $\mu_{20}^{(0)}\approx 12.98$.}
\label{fig:exploc_ellip}
\end{figure}

\begin{figure}[!ht]
\centering
\subfloat[\label{fig:coupea}]{\includegraphics[width=0.45\linewidth]{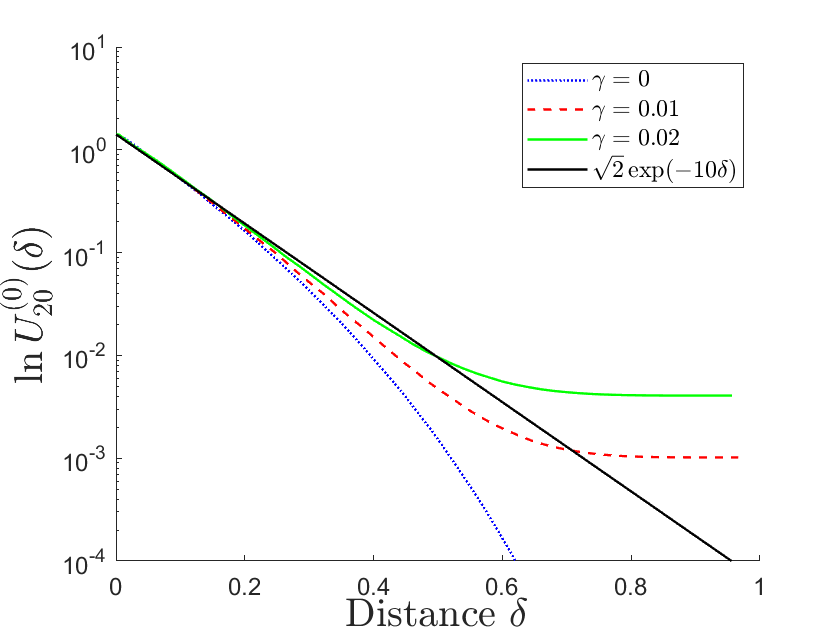}} \hskip1mm
\subfloat[\label{fig:coupeb}]{\includegraphics[width=0.45\linewidth]{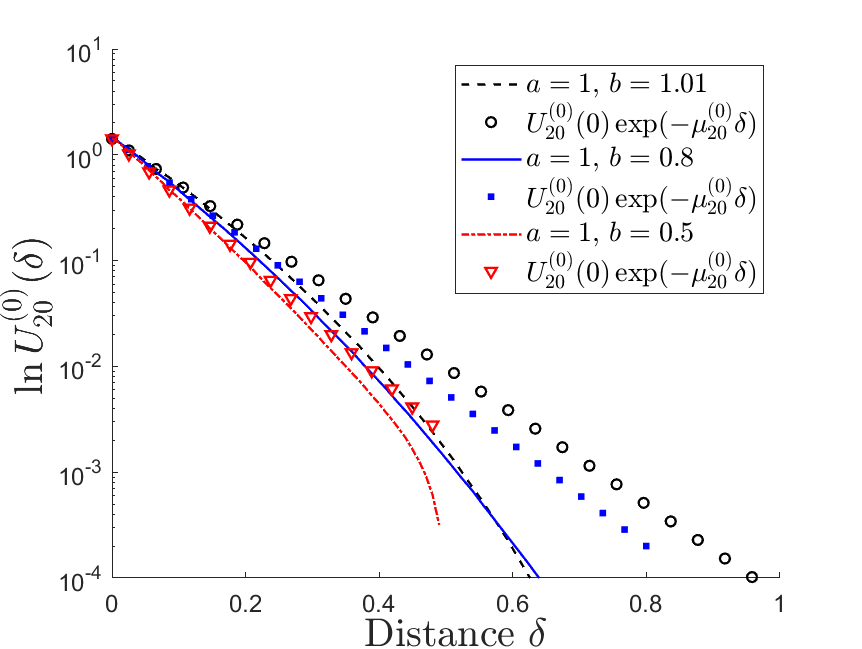}}
\caption{Log-plot of the decay of $U_{20}^{(0)}(\delta)$ for (a) a disk ($\gamma=0$) and two deformed disks ($\gamma=0.01$ and $\gamma=0.02$), where solid black line indicates the known decay $\sqrt{2}\exp(-10\delta)$ for a disk; (b) an ellipse with semiaxes $a=1$, $b=1.01$, an ellipse with semiaxes $a=1$, $b=0.8$, and an ellipse with semiaxes $a=1$, $b=0.5$, where lines indicate the expected decay $U_{20}^{(0)}(0) \exp(-\mu_{20}^{(0)} \delta)$. The related eigenvalues are $\mu_{20}^{(0)}\approx 9.95$,  $\mu_{20}^{(0)}\approx 11.08$ and $\mu_{20}^{(0)}\approx 12.98$ respectively.}
\label{fig:coupe}
\end{figure}

In order to investigate the exponential decay of eigenfunctions in the whole domain, we introduce the function
\begin{equation}
U_k^{(p)}(\delta) = \sqrt{|\pa|} \max\limits_{\x\in \gamma_{\delta}} |V_k^{(p)}(\x)|,
\end{equation}
where $\gamma_{\delta}$ is the contour line of points in $\Omega$ at distance $\delta$ from the boundary $\pa$.  
We aim to test whether the following approximation holds:
\begin{equation}\label{eq:testU}
U_k^{(p)}(\delta) \approx U_k^{(p)}(0) e^{-\mu_k^{(p)}\delta}.
\end{equation}
\Cref{fig:coupea} shows the log-plot of $U_{20}^{(0)}(\delta)$ for the three eigenfunctions shown on \cref{fig:exploc_smooth} and compare its exponential decay with that known for the disk. One sees that the more the disk is perturbed, the earlier the exponential decay $\exp(-\mu_k^{(0)}\delta)$ stops to approximate the behavior of the related eigenfunction $V_k^{(0)}$. This indicates that the determination of an optimal $\eta$ in \cref{eq:polt} might not be straightforward even for domains with real-analytic boundary and requires further analysis.
Note that the numerical results shown in \cref{fig:coupe} are more accurately described by the truncated exponential $\exp\left(-\mu_k^{(0)} \min\{\epsilon, |\x-\pa|\}\right)$ that appears in \cref{eq:helffer}. Indeed, there is a cutoff distance $\epsilon$ above which the exponential decay with the decay rate $\mu_k^{(0)}$ does not work. In other words, $\epsilon$ determines an inner ``central" region of the domain in which the eigenfunction decays slower. In turn, \cref{fig:coupeb} illustrates that there is no such a central region for the considered ellipses, i.e., there may not be need in the truncated exponential for these domains.  From a practical point of view, the natural question is to know whether the exponential function $b_k \exp(-\mu_k^{(0)} |\x-\pa|)$, or its truncated form $b_k \exp\left(-\mu_k^{(0)} \min\{\epsilon, |\x-\pa|\}\right)$, can be an accurate approximation of the Steklov eigenfunction $V_k^{(0)}$. For instance, can one choose $b_k$ such that $ \max\limits_{x\in\Omega}\left||V_k^{(p)}(\x)| -b_k\exp(-\mu_k^{(0)} |\x-\pa|)\right|$ is small enough?

Finally, we complete this section by providing a complementary insight onto the localization of $V_k^{(p)}$ at large $p$.
Multiplying \cref{prob:VkSteklov} by $V_k^{(p)}$, integrating over $\Omega$, using
the Green's formula, and employing the Steklov boundary condition, one
easily gets
\begin{equation}  \label{eq:VVnorm}
p\int_{\Omega} [V_k^{(p)}(\x)]^2 d\x + \int_{\Omega} |\nabla V_k^{(p)}(\x)|^2 d\x = \mu_k^{(p)} \,.
\end{equation}
The eigenvalue determines therefore the combination of
$L^2(\Omega)$-norms of $V_k^{(p)}$ and of its gradient.  It is
instructive to compute both norms separately.  In Appendix \ref{secap:norm}, we
derive the following relation
\begin{equation}   \label{eq:Vnorm}
\int_{\Omega} |V_k^{(p)}(\x)|^2 d\x = \partial_p \mu_k^{(p)} .
\end{equation}
Combining both equations, we also find
\begin{equation}
\int_{\Omega} |\nabla V_k^{(p)}(\x)|^2d\x  = \mu_k^{(p)} - p \, \partial_p \mu_k^{(p)} \,.
\end{equation}
As we discussed in \cref{sec:asymp}, the eigenvalues $\mu_k^{(p)}$ grow as
$\sqrt{p}$ at large $p$, and so does the left-hand side of
Eq. (\ref{eq:VVnorm}).  In turn, Eq. (\ref{eq:Vnorm}) implies that the
$L^2(\Omega)$-norm of $V_k^{(p)}$ vanishes as $p^{-1/4}$ as $p\to
\infty$.  This is a weaker form of the localization of the
Steklov eigenfunctions near the boundary.

\section{Discussion and conclusion}\label{sec:discussion}
In this paper, we numerically investigated the spectral properties of the Dirichlet-to-Neumann operator $\Mp$ and their dependence on the parameter $p$ and the domain geometry. We considered various shapes, including ellipses, triangles, rectangles, regular polygons and prefractal Koch snowflakes. 

Our first contribution concerned the asymptotic behavior of the eigenvalues of $\Mp$.
For all considered shapes, we confirmed the validity of the asymptotic relations \cref{eq:asymp-mu2,eq:asymp-ck} for large and small $p$. In the limit $p\to\infty$, the coefficients $c_k$ in \cref{eq:asymp-ck} were known to be 1 for bounded domains with $\mathcal{C}^1$ boundary. In turn, their values for polygonal domains were unknown. We conjectured that as $p$ increases, first Steklov eigenfunctions are getting localized near the corners of a polygonal domain, and suggested an iterative procedure to obtain the coefficients $c_k$. This conjectural relation was numerically validated on several domains. We are unaware of earlier studies on the coefficients $c_k$ for the asymptotic behavior of the eigenvalues $\mu_k^{(p)}$ of the Dirichlet-to-Neumann operator. However, additional insights can be gained from the related problem of the Robin Laplacian defined in \cref{prob:lap}.
In fact, Lacey \textit{et al.} considered the asymptotic behavior of the smallest eigenvalue $\lambda_0^{(q)}$ of the Robin Laplacian in the limit $q\to -\infty$, and found that \cite{lacey1998multidimensional}:
\begin{equation}\label{eq:lacey}
\lambda_0^{(q)} \simeq - q^2/\sin^2(\alpha/2)
\end{equation}
in a corner with angle $\alpha < \pi$.
This result was generalized by Levitin and Parnovski \cite{levitin2008principal} who proved \cref{eq:lacey} for a polygonal domain with angles
$\alpha_1, … \alpha_n$, such that $0 < \alpha_i < \pi$, and $\alpha = \min\{\alpha_i\}$.
Then, Khalile \cite{khalile2018spectral} extended this result for the first $n$ eigenvalues of
the Robin Laplacian for polygons with angles between $\pi/3$ and $\pi$.
Finally, Khalile and Pankrashkin \cite{khalile2018eigenvalues} considered a similar problem for an
infinite sector and showed the asymptotic behavior $\lambda_n^{(q)} \simeq  -
q^2/((2n+1)\alpha/2)^2$ in the limit $\alpha\to 0$.  The duality between the Robin problem
and the Steklov problem allows one to invert \cref{eq:lacey} to get $\mu_0^{(p)} \simeq \sin(\alpha/2) \sqrt{p}$ as $p\to\infty$ for a polygon. This asymptotic result, which is a direct consequence of \cite{levitin2008principal}, confirms our conjecture for the smallest eigenvalue $\mu_0^{(p)}$. In turn, our conjectural iterative procedure relates the asymptotic behavior of all eigenvalues to the angles of a polygonal domain. Its rigorous demonstration presents an intersting open problem. Perhaps, the most challeging part is to prove that our procedure yields all coefficients $c_k<1$.

The second result concerned some spectral expansions that appear in the theory of diffusion-controlled reactions \cite{grebenkov2020paradigm}. For symmetric domains such as ellipses, rectangles and even prefactal Koch snowflakes, we numerically observed that many coefficients $A_k^{(p)}$ vanish for a wide range of $k$ and $p$. We argued that this behavior is a consequence of domain symmetries. In turn, this property breaks for generic domains such as an arbitrary triangle. 
On one hand, the cancellation of many coefficients $A_k^{(p)}$ due to domain symmetries can considerably  simplify spectral expansions and the analysis of diffusion-controlled reactions in such domains. Moreover, these symmetries can potentially be used to design domains with specific properties. On the other hand, the use of symmetric domains as examples in theoretical and numerical studies may lead to erroneous conjectures and conclusions. For instance, the cancellation of coefficients $A_k^{(p)}$ with $k>0$ for the disk is a very specific consequence of its rotational symmetry, which fails for less symmetric domains. In particular, many coefficients $A_k^{(p)}$ contribute to the spectral expansion in the case of a generic triangle. It is therefore an open question how many terms are relevant and how their number depends on domain shape. More formally, it is known from the general arguments that $A_k^{(p)} \to 0$ as $k\to \infty$, but the speed of decay remains unknown. This question is practically important because a rapid decay of $A_k^{(p)}$ may allow truncating some spectral expansions to get useful approximations with a limited number of contributing eigenmodes.

Last, we investigated the localization of Steklov eigenfunctions $V_k^{(p)}$ on the boundary in the presence of corners. Even if the boundary is polygonal and thus is not real-analytic, \cref{fig:exploc_poly} illustrated the exponential decay of $|V_k^{(0)}|$ away from the boundary for a broad range of indices $k$. These numerical examples could motivate further mathematical investigations in this direction. However, we also observed that there exist domains and indices $k$, for which deviations from the exponential decay, characterized by $B_k(\x)$ and $U_k^{(p)}(\delta)$, become large, suggesting that the hypothesis $\eta\approx1$ in \cref{eq:polt} might not be relevant for all $k$ and in the whole domain. In particular, we noticed that only a slight perturbation of the disk can result in high values of $B_k(\x)$ located in the center of the domain, indicating that the related eigenfunctions decrease slower inside the domain than near the boundary. A more systematic study of this behavior through the estimation of the constants $B$ and $\eta$ in \cref{eq:polt} and their relation with the domain geometry present an interesting perspective to this work.

\section*{Acknowledgments}
The authors thank Prof. Iosif Polterovich, Prof. Michael Levitin and Prof. Bernard Helffer for fruitful discussions and suggestions.
The authors acknowledge Antoine Moutal for his early contributions in the numerical implementation of the finite element method.


\section{Appendix}

\subsection{Numerical validation: solution for rectangles}\label{secap:num}

Laugesen studied the spectral properties of the Robin Laplacian in
various rectangular domains (or cuboids) in $\mathbb{R}^d$ \cite{laugesen2019robin}. In this section, we adapt his analysis to the Steklov problem and give the explicit formulas for rectangles $\Omega = (0,b_1)\times (0,b_2)\subset \mathbb{R}^2$.
The separation of variables in the modified Helmholtz equation $(p-
\Delta) u = 0$ yields:
\begin{equation}  \label{eq:u_cuboid}
u(\x) = \prod\limits_{n=1}^2 \biggl(\alpha_n \cosh(\alpha_n x_n/b_n) - \mu b_n \sinh(\alpha_n x_n/b_n)\biggr),
\end{equation}
with
\begin{align}  
\tanh(\alpha_1) &= \frac{2\mu \alpha_1 b_1}{\alpha_1^2 + \mu^2 b_1^2}, \label{eq:alpha_i1} \\ 
\tanh(\alpha_2) &= \frac{2\mu \alpha_2 b_2}{\alpha_2^2 + \mu^2 b_2^2},  \label{eq:alpha_i2}
\end{align}
and
\begin{equation} \label{eq:alpha1_d}
\frac{\alpha_1^2}{b_1^2} + \frac{\alpha_2^2}{b_2^2} = p.
\end{equation}
Solving the system of three nonlinear equations (\ref{eq:alpha_i1}), (\ref{eq:alpha_i2}), (\ref{eq:alpha1_d}), one can determine the unknown coefficients
$\alpha_1,\alpha_2$ and the eigenvalue $\mu$. We
stress that solutions $\alpha_n$ can be either real, or purely
imaginary. In order to solve the system of three nonlinear equations, we consider \cref{eq:alpha_i1} as the quadratic equation on $\mu$,
whose two solutions are
\begin{align}  \label{eq:mu_pm}
\mu_+ & = \frac{\alpha_1}{b_1}\, \ctanh(\alpha_1/2), \quad
\mu_-  = \frac{\alpha_1}{b_1}\, \tanh(\alpha_1/2).
\end{align}

Let us first focus on the case $p = 0$, for which
\cref{eq:alpha1_d} implies $\alpha_2 = i \alpha_1 b_2/b_1$.
Substituting this expression into \cref{eq:alpha_i2}, we
get
\begin{equation}\label{eq:alpha1}
\tan(\alpha_1 b_2/b_1) = \frac{2\alpha_1 b_1\mu}{\mu^2b_1^2 - \alpha_1^2} \,.
\end{equation}
Substituting $\mu_\pm$ from Eq. (\ref{eq:mu_pm}) into this relation,
one has

\begin{equation}  \label{eq:alpha_rectangle0}
\tan(\alpha_1 b_2/b_1) = \pm \sinh(\alpha_1) .
\end{equation}
Each of these equations has infinitely many real solutions that
determine the eigenvalues of the Dirichlet-to-Neumann operator $\mathcal{M}_0$
for a rectangle.

When $b_1 = b_2 = 2$, one can set $\kappa = \alpha_1/2$ and use
trigonometric relations to rewrite Eq. (\ref{eq:alpha_rectangle0})
with plus sign in terms of $\kappa$ as:
\begin{equation*}
\frac{\tan(\kappa)}{1 - \tan^2(\kappa)} = \frac{\tanh(\kappa)}{1-\tanh^2(\kappa)} = \frac{-\ctanh(\kappa)}{1-\ctanh^2(\kappa)} \,,
\end{equation*}
which is equivalent to $\tan(\kappa) = \tanh(\kappa)$ or $\tan(\kappa)
= - \ctanh(\kappa)$.  Similarly, \cref{eq:alpha_rectangle0} with
minus sign reads
\begin{equation*}
\frac{\tan(\kappa)}{1 - \tan^2(\kappa)} = \frac{\ctanh(\kappa)}{1-\ctanh^2(\kappa)} = \frac{- \tanh(\kappa)}{1-\tanh^2(\kappa)} \,,
\end{equation*}
which is equivalent to $\tan(\kappa) = \ctanh(\kappa)$ or
$\tan(\kappa) = -\tanh(\kappa)$.  We retrieve therefore the equations
reported in Table 7.1 of \cite{levitin2023topics} for the square $(-1,1)^2$.  In
this case, there is an additional eigenvalue $1$, which corresponds to
the eigenfunction $x_1x_2$.  This eigenfunction is not included in the
general form (\ref{eq:u_cuboid}). 

\vskip 2mm
Now we return to the case $p > 0$.  We search for real solutions
$\alpha_1$.  Setting $\alpha_1 = (b_1/b_2) \sqrt{-\alpha_2^2 +p
b_2^2}$ and substituting $\mu_+$ from Eq. (\ref{eq:mu_pm}), we get the
following equation on $\alpha_2$:
\begin{equation}  \label{eq:tanh_alpha2}
\tanh(\alpha_2) = - \alpha_2 \frac{(\alpha_1 b_2/b_1) \sinh(\alpha_1)}{\alpha_2^2 -p b_2^2 \cosh^2(\alpha_1/2)} \,.
\end{equation}
There is a finite number of real solutions $\alpha_2$ of this equation
on the interval from $0$ to $b_2\sqrt{-\lambda}$, for which $\alpha_1$
is real.  In turn, there are infinitely many purely imaginary
solutions $\alpha_2$.  Setting $\alpha_2 = i\alpha$ and thus $\alpha_1
= (b_1/b_2) \sqrt{\alpha^2 +p b_2^2}$, we transform the above
equation into
\begin{equation}  \label{eq:auxil11}
\tan(\alpha) = \alpha \frac{(\alpha_1 b_2/b_1) \sinh(\alpha_1)}{\alpha^2 +p b_2^2 \cosh^2(\alpha_1/2)} \,.
\end{equation}
As $p> 0$, $\alpha_1$ is real and the right-hand side of
Eq. (\ref{eq:auxil11}) is positive.  This equation has infinitely many
solutions, which lie on the intervals $(\pi k,\pi k + \pi/2)$, with $k
= 0,1,2,\ldots$.

Similarly, using $\mu_-$ from Eq. (\ref{eq:mu_pm}), we get another
equation on $\alpha_2$:
\begin{equation}
\tanh(\alpha_2) = \alpha_2 \frac{(\alpha_1 b_2/b_1) \sinh(\alpha_1)}{\alpha^2 +p b_2^2 \sinh^2(\alpha_1/2)} \,,
\end{equation}
which has a finite number of real solutions $\alpha_2$ on the interval
$(0,b_2\sqrt{p})$.  In turn, setting again $\alpha_2 =
i\alpha$, one gets
\begin{equation*}
\tan(\alpha) = - \alpha \frac{(\alpha_1 b_2/b_1) \sinh(\alpha_1)}{\alpha^2 -p b_2^2 \sinh^2(\alpha_1/2)} \,.
\end{equation*}
Since the denominator in the right-hand side can change sign, it is
more convenient to rewrite this relation as
\begin{equation}
\frac{\alpha}{\tan(\alpha)} = - \frac{\alpha^2 -p b_2^2 \sinh^2(\alpha_1/2)}{(\alpha_1 b_2/b_1) \sinh(\alpha_1)} \,.
\end{equation}
There are infinitely many solutions of this equation.  

The above computation allows one to find all solutions, for which
$\alpha_1$ is real, while $\alpha_2$ is either real, or purely
imaginary.  Exchanging the roles of $\alpha_1$ and $\alpha_2$, one can
also determine the missing pairs, for which $\alpha_2$ is real while
$\alpha_1$ is either real, or purely imaginary.  Combining all these
solutions, we determine the eigenvalues of the Dirichlet-to-Neumann
operator $\Mp$ for the rectangle. 

\Cref{tab:recteigv} summarizes the first $11$ eigenvalues of $\Mp$ for a rectangles of sides $1$ and $2$, with $p=1$ (\cref{sfig:rect} shows the associated eigenfunction $V_4^{(1)}$). One sees that the eigenvalues in the third column, which were numerically obtained by Method 1, are in excellent agreement with the exact ones.

\begin{table}[!h]
\centering
$\begin{array}[!ht]{|c|c|c|c|}
\hline \textbf{Index k}  & \textbf{Exact} & \textbf {Method 1}  & \textbf {Method 2} \\
\hline 
0 & 0.3105 & 0.3105 &  0.3105 \\
1 & 0.7511 & 0.7511 &  0.7512\\
2 & 1.6451 & 1.6451 &  1.6435\\
3 & 1.7342 & 1.7342 &  1.7332\\
4 & 2.2304 & 2.2305 &  2.2315\\
5 & 2.9051 & 2.9051 &  2.9219\\
6 & 3.9156 & 3.9159 &  3.9232\\
7 & 4.1665 & 4.1668 &  4.2058\\
8 & 4.7951 & 4.7955 &  4.8040\\
9 & 4.7961 & 4.7966 &  4.8147\\
10 & 5.5419 & 5.5426 & 5.5582\\
\hline
\end{array}$
\caption{List of the first 11 eigenvalues $\mu_k^{(p)}$ for a rectangle with sides $1$ and $2$, with $p=1$. For Method 1, the mesh is composed of 179918 triangles, and the maximal mesh size is 0.005. For Method 2, the mesh is composed of 16256 triangles, the maximum mesh size is 0.03, and the series in \cref{eq:34} was truncated to 88 eigenfunctions of the Laplace operator.
}
\label{tab:recteigv}
\end{table}


Let us inspect the limit $p \to \infty$.  Setting $0 \leq
\gamma \leq \pi/2$, it is convenient to write $\alpha_1 = b_1 \sqrt{p}
\sin\gamma$ and $\alpha_2 = b_2 \sqrt{p} \cos\gamma$ that satisfy
Eq. (\ref{eq:alpha1_d}).  Substituting these expressions into
Eq. (\ref{eq:tanh_alpha2}), we get
\begin{equation}
\tanh(b_2 \sqrt{p}\cos\gamma) = \frac{\sin(2\gamma) \sinh(b_1\sqrt{p}\sin\gamma)}{2[\cosh^2(b_1\sqrt{p}\sin\gamma/2) - \cos^2(\gamma)]} \,.
\end{equation}
One can easily check that $\gamma \approx 0$ is not compatible with
this equation.  As a consequence, since $b_1 \sqrt{p} \sin\gamma \gg
1$ in the limit $p\to \infty$, the right-hand side of the above
equation is close to $\sin(2\gamma)$.  If $b_2 \sqrt{p} \cos \gamma
\gg 1$, one gets $1 \approx \sin(2\gamma)$, from which $\gamma =
\pi/4$ and thus $\mu \approx \alpha_1/b_1 \approx
\sin(\pi/4)\sqrt{p}$.  In turn, if $b_2 \sqrt{p} \cos \gamma
\ll 1$, then $\gamma \approx \pi/2$, and thus $\mu \approx \sqrt{p}$.
These qualitative arguments are consistent with our numerical
predictions on the behavior of the eigenvalues in the limit $p\to
\infty$.  However, more accurate analysis is needed to claim that only
two eigenvalues behave as $\mu \approx \sin(\pi/4)\sqrt{p}$.

\subsection{Localization in rectangles}
\label{secap:loc-rect}
Let us focus on $p = 0$ and consider the solutions of \cref{eq:alpha1} for which
$\alpha_1$ is real, while $\alpha_2 = i\alpha_1 b_2/b_1$ is purely
imaginary. The associated eigenfunction, which is given by \cref{eq:u_cuboid} up to a normalization, is factored as $u_1(x_1) u_2(x_2)$.  As the factor
$u_2(x_2)$ exhibits oscillatory behavior along $x_2$ coordinate, we
focus on $u_1(x_1)$, which can be written as
\begin{equation}
u_1(x_1) = \frac{\alpha_1 - \mu b_1}{2} e^{\alpha_1 x_1/b_1} + \frac{\alpha_1 + \mu b_1}{2} e^{-\alpha_1 x_1/b_1} .
\end{equation}
Using the first relation in \cref{eq:mu_pm} to express $\mu$, one gets then
\begin{equation}\label{eq:um_rect}
u_1(x_1) = \frac{\alpha_1}{1-e^{-\alpha_1}} \biggl( e^{- \alpha_1 x_1/b_1} - e^{-\alpha_1 (1-x_1/b_1)}\biggr) ,
\end{equation}
i.e., $|u_1(x_1)|$ decays exponentially with the distance $\delta_1 =
\min\{x_1, b_1-x_1\}$ from either of two endpoints of the interval
$(0,b_1)$:
\begin{equation}\label{eq:upb1}
|u_1(x_1)| \leq \frac{\alpha_1}{1-e^{-\alpha_1}} e^{- \delta_1 \alpha_1/b_1} \simeq \frac{\alpha_1}{1-e^{-\alpha_1}} e^{-\mu_+ \delta_1},
\end{equation}
where $\alpha_1/b_1$ is exponentially close to $\mu_+$ when $\alpha_1$
is large enough.  

Similarly, if one uses the second relation in \cref{eq:mu_pm} to express
$\mu$, one gets
\begin{equation}\label{eq:up_rect}
u_1(x_1) = \frac{\alpha_1}{1+e^{-\alpha_1}} \biggl( e^{- \alpha_1 x_1/b_1} + e^{-\alpha_1 (1-x_1/b_1)}\biggr) ,
\end{equation}
which can bounded as
\begin{equation}\label{eq:upb2}
|u_1(x_1)| \leq \frac{2\alpha_1}{1+e^{-\alpha_1}} e^{- \delta_1 \alpha_1/b_1} \simeq \frac{2\alpha_1}{1+e^{-\alpha_1}} e^{-\mu_- \delta_1}.
\end{equation}

Setting $\delta_2 = \min\{x_2, b_2-x_2\}$, one gets $|\x-\pa| = \min\{\delta_1,\delta_2\}\leq\delta_1$, so that the upper bounds \cref{eq:upb1,eq:upb2} can be written as 
\begin{equation}
|u_1(x_1)| \leq C_1 e^{-\mu |\x-\pa|},
\end{equation}
with a constant $C_1$. Finally, the analytic function $u_2(x_2)$ is bounded by its maximum, one has 
\begin{equation}\label{eq:V_rect}
|\sqrt{\pa}||V_k^{(0)}| \lesssim B \exp\left(-\mu_k^{(0)}|\x-\pa|\right),
\end{equation}
i.e. we retrieve the upper bound (\ref{eq:polt}) with $\eta\approx1$. 

\begin{figure}[!ht]
\centering
\subfloat[\label{sfig:rect2_U10}]{\includegraphics[width=0.33\linewidth]{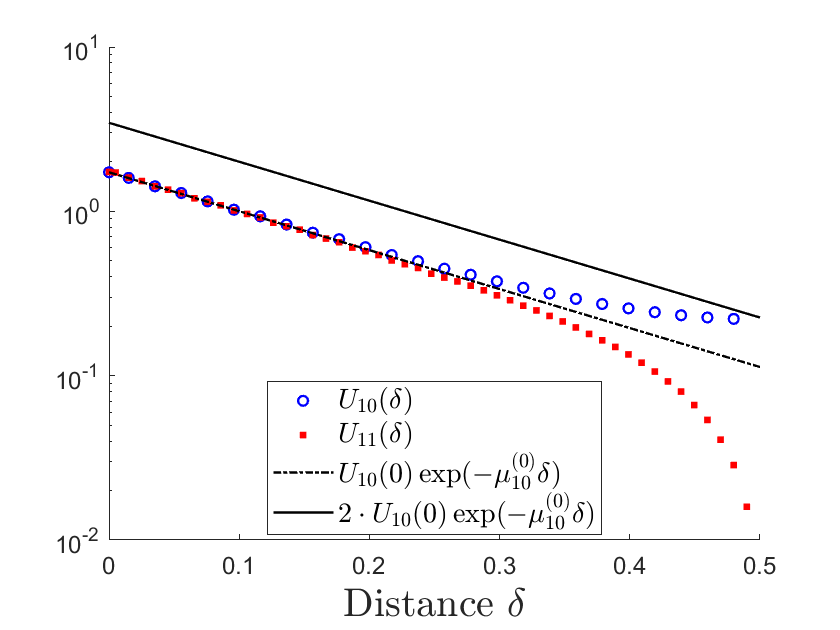}}
\subfloat[\label{sfig:rect2_V20}]{\includegraphics[width=0.33\linewidth]{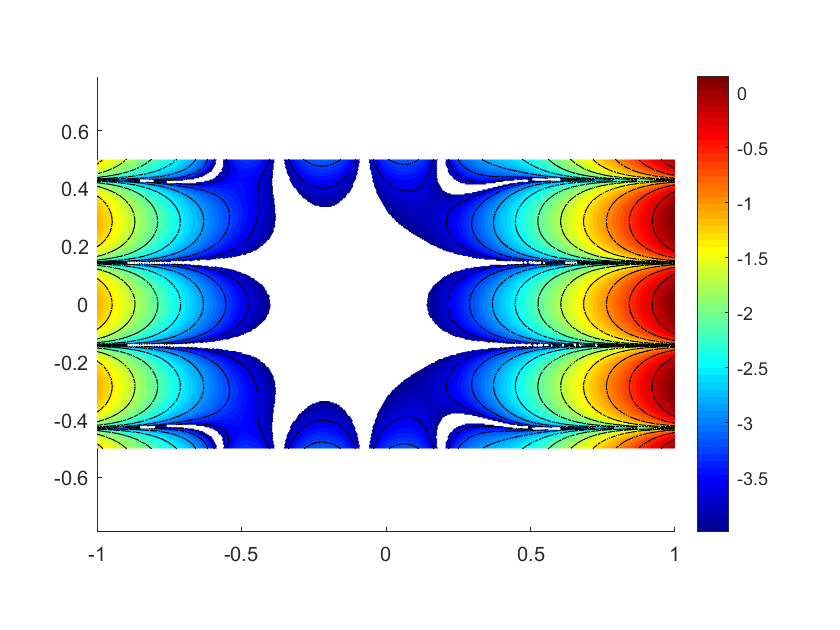}}
\subfloat[\label{sfig:rect2_V21}]{\includegraphics[width=0.33\linewidth]{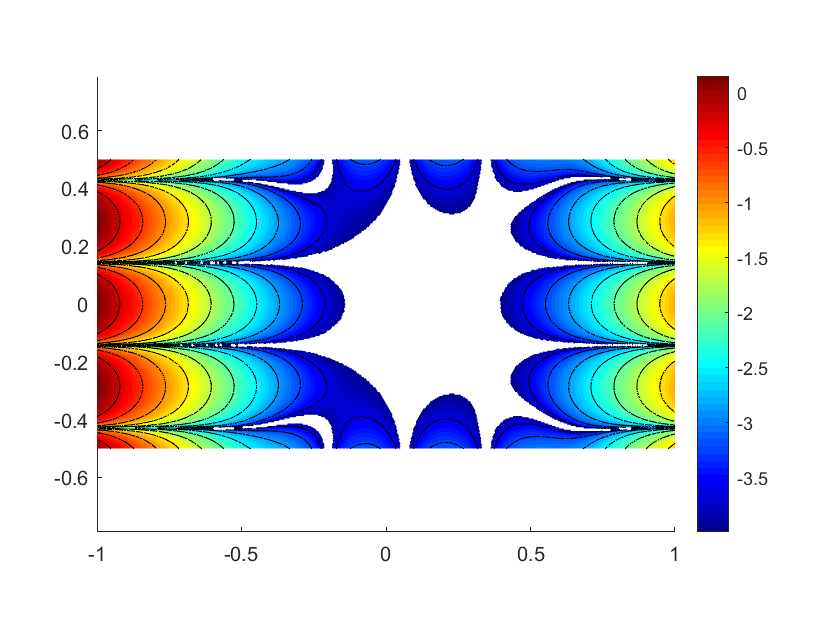}}
\caption{(a) Log-plot of $U_{10}^{(0)}(\delta)$ and $U_{11}^{(0)}(\delta)$ for a rectangle with sides $2$ and $1$, the related eigenvalues are $\mu_{10}^{(0)}\approx5.46$ and $\mu_{11}^{(0)}\approx5.54$, and lines indicate the expected decay $U_{10}^{(0)}(0) \exp(-\mu_{10}^{(0)} \delta)$ and the upper bound $2\cdot U_{10}^{(0)}(0) \exp(-\mu_{10}^{(0)} \delta)$ for one of the eigenfunctions; (b) and (c) the localization and decay of $V_{20}^{(0)}$ and $V_{21}^{(0)}$ away from the boundary. Colormaps illustrate the behavior of $\lg(|V_{20}^{(0)}|)$ on (b), and $\lg(|V_{21}^{(0)}|)$ on (c). The related eigenvalues are $\mu_{20}^{(0)} = \mu_{21}^{(0)}\approx11.00$.}
\label{fig:exploc_rect}
\end{figure}

We note that \cref{eq:um_rect,eq:up_rect} highlight the expected symmetry of the eigenfunctions discussed in \cref{sec:trunc}, namely, $u_1(x_1) = - u_1(b_1-x_1)$ in \cref{eq:um_rect} and $u_1(x_1) = u_1(b_1-x_1)$ in \cref{eq:up_rect}. In the former case, the antisymmetric eigenfunction vanishes in the middle so that it may decay even faster than exponential near the middle point. In turn, the symmetric eigenfunction does not vanish at the middle, and both terms in \cref{eq:up_rect} provide equal contributions in the middle. Such a function is expected to decay slower in the middle. \Cref{sfig:rect2_U10} illustrates this behavior for eigenfunctions $V_{10}^{(0)}$ and $V_{11}^{(0)}$.  
One sees that, even though the exponentially decaying upper bound \cref{eq:V_rect} holds for both cases, it does not necessarily approximate the eigenfunction. In fact, to fulfill the upper bound for symmetric eigenfunctions, we had to add the factor 2 in \cref{eq:upb2}, which shifts the upper bound from the expected exponential behavior of the eigenfunction. 
The above symmetry argument is applicable to eigenfunctions corresponding to simple (non-degenerate) eigenvalues. In turn, if two (or more) eigenfunctions correspond to the same eigenvalue, then their linear superposition is also an eigenfunction that can break this symmetry. We illustrate this situation for a rectangle with sides 2 and 1, for which \cref{sfig:rect2_V20,sfig:rect2_V21} shows two eigenfunctions $V_{20}^{(0)}$ and $V_{21}^{(0)}$ that correspond to  a twice degenerate eigenvalue. One of these eigenfunctions is localized on the left edge and the other is on the right edge of the rectangle. 

Refining the above arguments, one can achieve more rigorous statements, in particular on the closeness of $\eta$ to 1. We expect that similar analysis can be performed for $p>0$. 

\subsection{Derivation of the $L^2(\Omega)$-norm of $V_k^{(p)}$}
\label{secap:norm}
In this Appendix, we derive the identity
(\ref{eq:Vnorm}).  First, applying the spectral expansion \cref{eq:34} and the
orthogonality of Laplacian eigenfunctions $u_k^{(q)}$ to each other,
we get a simple identity for Green's functions
\begin{equation}
\int_{\Omega} \tilde{G}_q(\x_1,p|\x) \tilde{G}_q(\x_2,p|\x) d\x = - \partial_p \tilde{G}_q(\x_2,p|\x_1)  \qquad (\forall ~ \x_1,\x_2\in\Omega).
\end{equation}
Next, using the representation (\ref{eq:Vk_G0}), one can write the squared
$L^2(\Omega)$-norm of the Steklov eigenfunction $V_k^{(p)}$ as
\begin{align*}
\| V_k^{(p)}\|_{L^2(\Omega)}^2 & = \int_{\Omega} |V_k^{(p)}(\x)|^2 d\x 
=  \int_{\pa}\mu_k^{(p)} \, v_k^{(p)}(\x_1)  d\x_1  \int_{\pa} \mu_k^{(p)} \, [v_k^{(p)}(\x_2)]^* d\x_2 \\
& \times \underbrace{\int_{\Omega} \tilde{G}_0(\x_1,p|\x) \, \tilde{G}_0(\x_2,p|\x)d\x}_{= -\partial_p \tilde{G}_0(\x_2,p|\x_1)} \\
& =  - [\mu_k^{(p)}]^2 \int_{\pa} v_k^{(p)}(\x_1) d\x_1 \int_{\pa}  [v_k^{(p)}(\x_2)]^*
\partial_p \biggl(\sum\limits_{j=0}^\infty \frac{[v_j^{(p)}(\x_1)]^* \, v_j^{(p)}(\x_2)}{\mu_j^{(p)}}\biggr) d\x_2,
\end{align*}
where we substituted the expansion (\ref{eq:G0}).  Writing
\begin{equation*}
\partial_p \biggl(\frac{[v_j^{(p)}(\x_1)]^* \, v_j^{(p)}(\x_2)}{\mu_j^{(p)}}\biggr)
= \frac{\partial_p [v_j^{(p)}(\x_1)]^* \, v_j^{(p)}(\x_2)}{\mu_j^{(p)}} + \frac{[v_j^{(p)}(\x_1)]^* \, \partial_p v_j^{(p)}(\x_2)}{\mu_j^{(p)}} 
- \frac{[v_j^{(p)}(\x_1)]^* \, v_j^{(p)}(\x_2)}{[\mu_j^{(p)}]^2} \partial_p \mu_k^{(p)} ,
\end{equation*}
one can separately evaluate three contributions by using the
orthonormality of eigenfunctions $\{v_k^{(p)}\}$ on $\pa$:
\begin{alignat*}{3}
\| V_k^{(p)}\|_{L^2(\Omega)}^2  &= - [\mu_k^{(p)}]^2 \biggl\{ \biggl( \int_{\pa} v_k^{(p)}(\x_1) 
\frac{\partial_p [v_k^{(p)}(\x_1)]^*}{\mu_k^{(p)}}d\x_1\biggr) \\
+ \biggl(\int_{\pa}  [v_k^{(p)}(\x_2)]^*  
\frac{\partial_p v_k^{(p)}(\x_2)}{\mu_k^{(p)}} d\x_2\biggr) - \frac{\partial_p \mu_k^{(p)}}{[\mu_j^{(p)}]^2}\biggr\} \\
& = - \mu_k^{(p)} \underbrace{\biggl( \int_{\pa} \partial_p |v_k^{(p)}(\x)|^2 d\x\biggr)}_{=0} 
+ \partial_p \mu_k^{(p)} .
\end{alignat*}
Exchanging the order of integration over $\x$ and differentiation with
respect to $p$ yields $\partial_p \| v_k^{(p)}\|_{L^2(\pa)}^2 = 0$ and
thus implies the identity (\ref{eq:Vnorm}).
\bibliographystyle{elsarticle-num}
\bibliography{biblio}

\end{document}